\documentclass{amsart}
\usepackage{amsmath}
\usepackage{amssymb} 
\usepackage{amscd} 
\usepackage{graphicx} 
\usepackage{epsfig}
\usepackage[dvips]{color}
\usepackage{slashbox}
\usepackage{bm}

\newtheorem{theorem}{Theorem}[section]
\newtheorem{lemma}[theorem]{Lemma}

\newtheorem{proposition}[theorem]{Proposition}
\newtheorem{corollary}[theorem]{Corollary}
\newtheorem{claim}[theorem]{Claim}

\theoremstyle{definition}
\newtheorem{definition}[theorem]{Definition}
\newtheorem{example}[theorem]{Example}

\theoremstyle{remark}
\newtheorem{remark}[theorem]{Remark}

\newtheorem{question}[theorem]{Question}

\numberwithin{equation}{section}
\numberwithin{figure}{section}
\numberwithin{table}{section}

\newcommand{\QED}[1]{\hspace*{\fill} $\square$(#1)\newline}

\begin{document}

\title[Seiferters for a trefoil knot]
{Neighbors of Seifert surgeries on a trefoil knot in the Seifert Surgery Network}

\author{Arnaud Deruelle, 
Katura Miyazaki,
and Kimihiko Motegi \\
}

\address{Institute of Natural Sciences, 
Nihon University, 
Tokyo 156-8550, Japan}
\email{aderuelle@math.chs.nihon-u.ac.jp}
\address{Faculty of Engineering, Tokyo Denki University, Tokyo 120-8551, 
Japan}
\email{miyazaki@cck.dendai.ac.jp}
\address{Department of Mathematics, Nihon University, 
Tokyo 156-8550, Japan}
\email{motegi@math.chs.nihon-u.ac.jp}

\date{}
\dedicatory{Dedicated to Fico Gonz\'alez Acu\~na on the occasion of his 70th birthday}

\begin{abstract}
A Seifert surgery is a pair $(K, m)$ of a knot $K$ in $S^3$
and an integer $m$ such that 
$m$--Dehn surgery on $K$ results in
a Seifert fiber space
allowed to contain fibers of index zero.
Twisting $K$ along a trivial knot called a seiferter for $(K, m)$
yields Seifert surgeries.
We study Seifert surgeries obtained from
those on a trefoil knot by twisting along their seiferters.
Although Seifert surgeries on a trefoil knot are the most basic ones,
this family is rich in variety.
For any $m \ne -2$ it contains
a successive triple of Seifert surgeries
$(K, m)$, $(K, m +1)$, $(K, m +2)$
on a hyperbolic knot $K$,
e.g.\ $17$--, $18$--, $19$--surgeries on
the $(-2, 3, 7)$ pretzel knot.
It contains infinitely many Seifert surgeries on strongly invertible
hyperbolic knots none of which arises from
the primitive/Seifert-fibered construction,
e.g.\ $(-1)$--surgery on the $(3, -3, -3)$ pretzel knot. 
\end{abstract}

\maketitle
{
\renewcommand{\thefootnote}{}
}

\section{Introduction}
\label{section:Introduction}
A pair $(K, m)$ of a knot $K$ in $S^3$ and an integer $m$ is
a \textit{Seifert surgery} if the result $K(m)$ of $m$--Dehn surgery
is a Seifert fiber space
which may contain a fiber of index $0$,
i.e.\ a degenerate fiber.
In this paper we allow Seifert fibrations to contain
degenerate fibers.
If $K(m)$ admits a degenerate Seifert fibration,
it is either a lens space or a connected sum of two lens spaces
\cite[Proposition~2.8(2), (3)]{DMM1}.
For a Seifert surgery $(K, m)$, when $K(m)$ admits
a non--degenerate Seifert fibration 
(i.e.\ $K(m)$ is not a connected sum of two lens spaces),
to emphasize this fact we also say that
$(K, m)$ is a \textit{Seifert fibered surgery}.

In \cite{DMM1}, we relate Seifert surgeries by
twists along ``seiferters" and define a $1$--dimensional complex
called the Seifert Surgery Network.
We briefly review the definition of the network.
Let $(K, m)$ be a Seifert surgery.
A knot $c \subset S^3 -K$ is a \textit{seiferter} for $(K, m)$
if $c$ is a trivial knot in $S^3$ but becomes a fiber in
a Seifert fibration of $K(m)$.
Let $K_p$ and $m_p$ be the images of $K$ and $m$
under twisting $p$ times along $c$;
in fact, $m_p = m +p(\mathrm{lk}(K, c))^2$.
Then, $(K_p, m_p)$ remains a Seifert surgery,
and the image of $c$ under the twisting is
a seiferter for $(K_p, m_p)$;
see the commutative diagram below.
We also consider twists along an ``annular pair of seiferters".
For two seiferters $c_1, c_2$ for $(K, m)$,
if $c_1$ and $c_2$ are fibers in the same Seifert fibration of $K(m)$,
then the (unordered) pair $\{c_1, c_2\}$ is
a \textit{pair of seiferters}.
A pair of seiferters is an \textit{annular pair of seiferters} if
$c_1$ and $c_2$ cobound an annulus $A$ in $S^3$.
After twisting along the annulus $A$
the images of $(K, m)$ and $\{c_1, c_2\}$
remain a Seifert surgery and an annular pair of seiferters for it.
The vertices of the \textit{Seifert Surgery Network}
 are all Seifert surgeries,
and two vertices of the network are connected by an edge if
one vertex (Seifert surgery) is obtained from the other
by a single twist along a seiferter or an annulus cobounded by
an annular pair of seiferters.
Refer to \cite[Subsection~2.4]{DMM1} for details of the definition.

\begin{eqnarray*}
\begin{CD}
	(K, m) @>\textrm{twist along } c\ (\textrm{resp.\ } \{c_1, c_2\})>> (K_p, m_p) \\
@V{m \textrm{--surgery on } K}VV
		@VV{m_p \textrm{--surgery on } K_p}V \\
	K(m) @>>\textrm{surgery along } c\ (\textrm{resp.\ } c_1 \cup c_2)> K_p(m_p)
\end{CD}
\end{eqnarray*}
\vskip 0.2cm
\begin{center}
\textsc{Diagram 1.} 
\end{center}

\begin{remark}
\label{rem:irrelevant}
\begin{enumerate}
\item
In \cite{DMM1}, an annular pair $\{c_1, c_2\}$ is defined
to be an ordered pair of $c_1$ and $c_2$ to specify
the direction of twist along the annulus
cobounded by $c_1 \cup c_2$. 
However, since we do not perform annulus twists in this paper,
annular pairs are presented as unordered pairs.
\item
If a seiferter $c$ for $(K, m)$ bounds a disk in $S^3 -K$,
we call $c$ \textit{irrelevant} and
do not regard it as a seiferter.
This is because no twists along irrelevant seiferters change
Seifert surgeries.
However, for pairs of seiferters $\{c_1, c_2\}$
we allow $c_i$ to be irrelevant.
Let $\{c_1, c_2\}$ be an annular pair for $(K, m)$.
If either $c_1$ and $c_2$ cobound an annulus disjoint from $K$ or
there is a 2--sphere in $S^3$ separating $c_i$ and $c_j\cup K$,
then twists along $\{c_1, c_2\}$ do not change $(K, m)$ or
have the same effect on $K$ as twists along $c_j$.
We thus call such an annular pair \textit{irrelevant},
and exclude it from annular pairs of seiferters.
\end{enumerate}
\end{remark}

\begin{figure}[htbp]
\begin{center}
\includegraphics[width=0.27\linewidth]{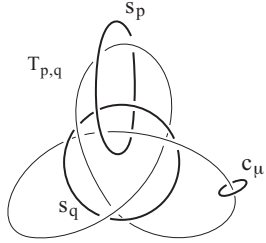}
\caption{Basic seiferters for $T_{p, q}$, where $p=-3, q=2$}
\label{basic}
\end{center}
\end{figure}

Any integral surgery on a torus knot $T_{p,q}$
$(|p| > q \ge 1)$ has at least three seiferters.
Let $s_p$ and $s_q$ be the exceptional fibers of indices $|p|$ and $q$
in the Seifert fibration of the exterior of $T_{p, q}$, respectively;
see Figure~\ref{basic}.
Since the Seifert fibration of the exterior extends to
$T_{p, q}(m)$ for any integer $m$,
the trivial knots $s_p$, $s_q$ are seiferters for $(T_{p,q}, m)$.
Furthermore, a meridian $c_{\mu}$ of $T_{p,q}$ is also a seiferter
for $T_{p, q}(m)$
because $c_{\mu}$ is isotopic to the core of the filled solid torus
in $T_{p, q}(m)$.
The seiferters $s_p$, $s_q$, $c_{\mu}$ are fibers 
of indices $|p|$, $q$, $|pq -m|$ in $T_{p, q}(m)$, respectively.
We call them \textit{basic seiferters} for $T_{p,q}$.
Note that $s_p$, $s_q$, $c_{\mu}$ in Figure~\ref{basic}
are fibers in a Seifert fibration of 
$T_{p, q}(m)$, simultaneously,
and any two of these seiferters cobound an annulus in $S^3$.
Thus, $\{ s_p, s_q\}$, $\{ s_p, c_{\mu} \}$, $\{ s_q, c_{\mu} \}$
in Figure~\ref{basic} are annular pairs of seiferters for
$(T_{p, q}, m)$, called \textit{basic annular pairs}.

In the network,
a path from $(T_{p, q}, m)$ to $(K, m')$
tells how the Seifert surgery $(K, m')$
is obtained from $(T_{p, q}, m)$
by a sequence of twistings along seiferters
and/or annular pairs of seiferters.
However, we cannot obtain a non-torus knot
by twisting a torus knot
along its basic seiferters or basic annular pairs.
To obtain a Seifert surgery on a hyperbolic knot
we need to twist along a ``hyperbolic seiferter".
A seiferter $c$ (resp.\ an annular pair $\{c_1, c_2\}$) for $(K, m)$ is 
\textit{hyperbolic} if $S^3 - K \cup c$
(resp.\ $S^3 - K \cup c_1 \cup c_2$) admits
a complete hyperbolic metric of finite volume. 
Twists along a ``hyperbolic seiferter"
or a ``hyperbolic annular pair" yield infinitely many
Seifert surgeries on hyperbolic knots.
We denote by $\mathcal{N}(T_{-3, 2})$
the set of Seifert surgeries obtained from
$(T_{-3, 2}, m)$ $(m \in \mathbb{Z})$
by twisting arbitrary times along seiferters or
annular pairs for $(T_{-3, 2}, m)$.

In this paper,
we find hyperbolic seiferters and hyperbolic annular pairs
of seiferters for $(T_{-3, 2}, m)$,
and study Seifert surgeries on hyperbolic knots 
that belong to $\mathcal{N}(T_{-3, 2})$.
We construct hyperbolic seiferters (resp.\ hyperbolic annular pairs)
for $(T_{-3, 2}, m)$ by applying ``$m$--moves" to
basic seiferters (resp.\ basic annular pairs) for $T_{-3,2}$.
An $m$--move is, in fact, a Kirby calculus handle-slide over
an $m$--framed knot,
and the definition is given in Section~\ref{section:m-move}.
Theorem~\ref{main} below follows from
Corollaries~\ref{hyperbolic seiferters},
\ref{hyperbolic annular pairs T-32}.

\begin{theorem}
\label{main}
A Seifert surgery $(T_{-3, 2}, m)$ has a hyperbolic seiferter for any integer $m \ne -4$; 
$(T_{-3, 2}, -4)$ has at least six hyperbolic annular pairs of seiferters.
Furthermore, 
if $m \le -8$ or $-1 \le m$, 
then $(T_{-3, 2}, m)$ has at least three hyperbolic seiferters and nine hyperbolic annular pairs of seiferters. 
\end{theorem}

The following two theorems are about surgeries
belonging to $\mathcal{N}(T_{-3, 2})$. 
A small Seifert fiber space is a 3--manifold which admits
a non-degenerate Seifert fibration over the $2$--sphere
containing exactly three exceptional fibers.
We call a Seifert surgery $(K, m)$
a \textit{small Seifert fibered surgery} if $K(m)$ is
a small Seifert fiber space.

\begin{theorem}
\label{maintriple}
For any integer $m$, there is a hyperbolic knot
whose $m$--, \mbox{$(m+1)$--~,} $(m+2)$--surgeries are
small Seifert fibered surgeries.
If $m\ne -2$, then such three successive surgeries can be 
found in $\mathcal{N}( T_{-3, 2} )$
$($Theorem~\ref{triple}$)$.
\end{theorem}

\begin{theorem}
\label{mainnonPS}
The neighborhood $\mathcal{N}(T_{-3, 2})$ contains infinitely many
small Seifert fibered surgeries on strongly invertible
hyperbolic knots 
which do not arise from the primitive/Seifert-fibered construction
introduced in \cite{D2}
$($Theorem~\ref{generalizedP(3,-3,-3)}$)$. 
\end{theorem}

\begin{figure}[htbp]
\begin{center}
\includegraphics[width=1.0\linewidth]{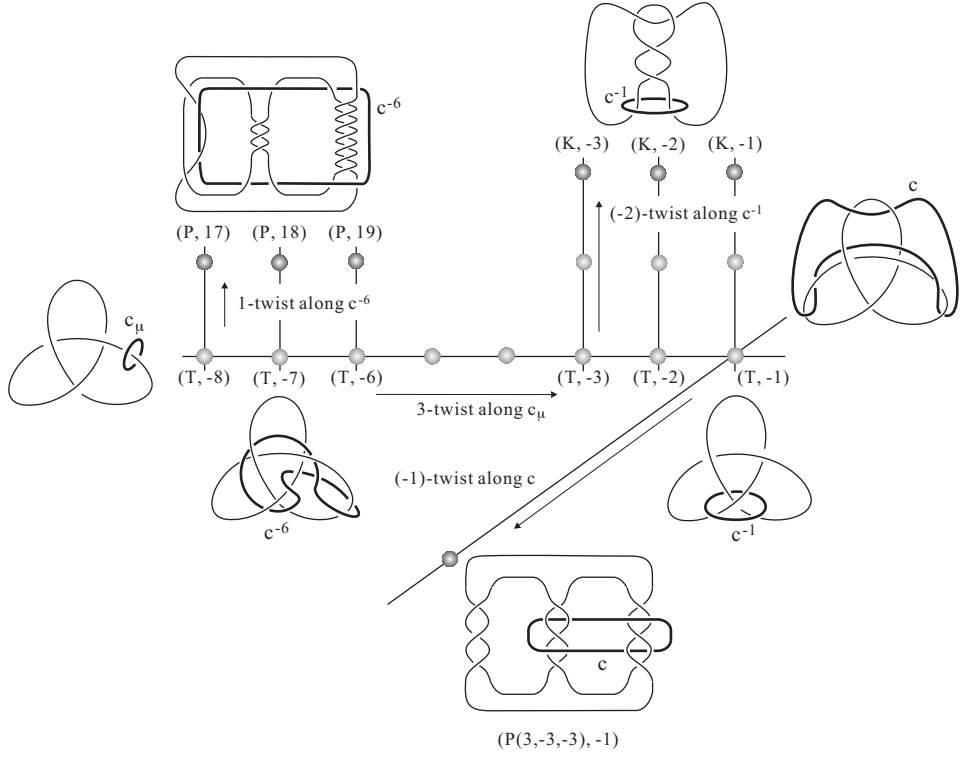}
\caption{$T = T_{-3, 2}$, $K$ is the figure-eight knot, and $P$
is the $(-2, 3, 7)$ pretzel knot.}
\label{figeightP}
\end{center}
\end{figure}

Figure~\ref{figeightP} is a portion of the subnetwork
$\mathcal{N}(T_{-3, 2})$.
Twists along the meridian $c_{\mu}$ generate the horizontal line in
Figure~\ref{figeightP}, 
which contains all integral surgeries on $T = T_{-3, 2}$. 
The trivial knots $c^m$ $(m =-1, -6)$ in Figure~\ref{figeightP}
are hyperbolic seiferters for
$(T_{-3, 2}, m)$, $(T_{-3, 2}, m-1)$, $(T_{-3, 2}, m-2)$
(Corollary~\ref{seiferter cm}, Proposition~\ref{cm is nonbasic});
$c$ in Figure~\ref{figeightP} is
a hyperbolic seiferter for $(T_{-3, 2}, -1)$ (Lemma~\ref{P(3,-3,-3)}).
Note that in Figure~\ref{figeightP} the images of
seiferters under twisting are denoted by the same symbols
as originals.
The $(-2)$--twist on $T_{-3, 2}$ along $c^{-1}$ yields 
the figure-eight knot $K$. 
Thus $(-2)$--twist along $c^{-1}$ converts $(T_{-3, 2}, m)$
to $(K, m -2w^2)$,
where $m =-1, -2, -3$ and $w =\mathrm{lk}(T_{-3, 2}, c^{-1}) =0$.
The $1$--twist along $c^{-6}$ yields
the $(-2, 3, 7)$ pretzel knot $P$,
so that the $1$--twist along $c^{-6}$
converts $(T_{-3, 2}, m)$ to $(P, m + w^2)$,
where $m =-6, -7, -8$ and $w = \mathrm{lk}(T_{-3, 2}, c^{-6}) =5$.
The slanted line through $(T_{-3, 2}, -1)$ is generated by
twists along the seiferter $c$.
The $1$--twist along $c$ yields $(-1)$--surgery on
the $(3, -3, -3)$ pretzel knot.
All Seifert surgeries on the slanted line
except $(T_{-3, 2}, -1)$ do not arise from
the primitive/Seifert-fibered construction.

\section{Preliminaries}
\label{section:m-move}

In this section we recall some results on $m$--moves to 
seiferters and annular pairs.

\begin{figure}[htbp]
\begin{center}
\includegraphics[width=1.0\linewidth]{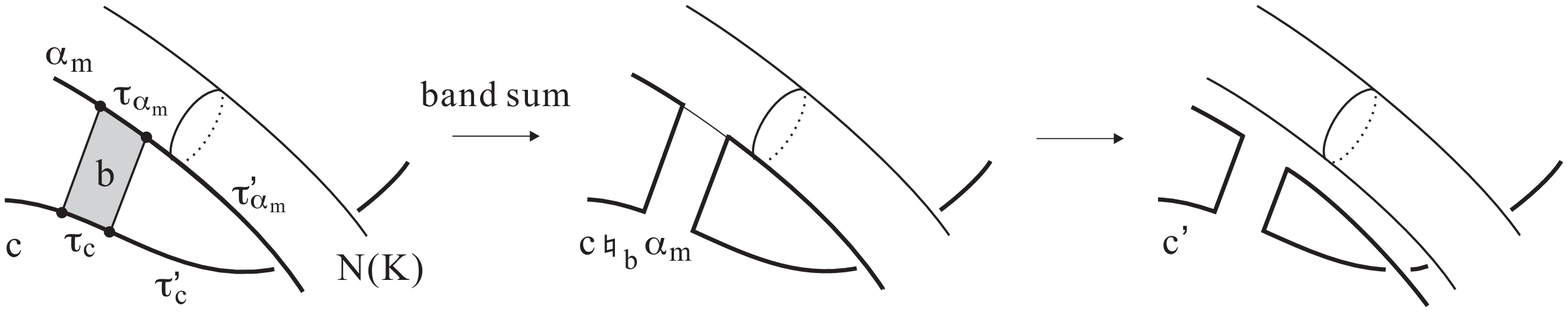}
\caption{$m$--move}
\label{movepush}
\end{center}
\end{figure}

\begin{definition}[$m$--moves]
Let $K$ be a knot in $S^3$ with a tubular neighborhood $N(K)$, 
and $c$ a knot in $S^3 - N(K)$. 
Take a simple closed curve $\alpha_m$ on $\partial N(K)$ 
representing a slope $m$. 
Let $b$ be a band in $S^3 -\mathrm{int}N(K)$
connecting $\alpha_m$ and $c$, and let
$b \cap \alpha_m = \tau_{\alpha_m}$, $b \cap c = \tau_c$.
We set $\tau'_{\alpha_m} = \alpha_m - \mathrm{int}\tau_{\alpha_m}$
and $\tau'_{c} = c - \mathrm{int}\tau_{c}$. 
Then the band connected sum 
$c\,\natural_b\,\alpha_m
= \tau'_c \cup (\partial b - \mathrm{int}(\tau_c \cup \tau_{\alpha_m}))
\cup \tau'_{\alpha_m}$ is a knot in $S^3 - \mathrm{int}N(K)$.  
Pushing $c\,\natural_b\,\alpha_m$ away 
from $\partial N(K)$, 
we obtain a knot $c'$ in $S^3 - N(K)$; 
see Figure~\ref{movepush}.
We say that $c'$ is obtained from $c$ by an \textit{$m$--move}
using the band $b$. 
\end{definition}

\begin{proposition}[{\cite[Propositions~2.19(3), 2.22]{DMM1}}]
\label{move}
Let $K$ be a nontrivial knot in $S^3$. 
Suppose that $(K, m)$ is a Seifert surgery with a seiferter $c$. 

\begin{enumerate}
\item
Assume that $c'$ is obtained from $c$ by a finite sequence of $m$--moves. 
Then, $c'$ is isotopic to $c$ in $K(m)$. 
Moreover, if $c'$ is unknotted in $S^3$, 
then $c'$ is also a seiferter for $(K, m)$. 
\item 
If $c'$ is obtained from $c$ by a single $m$--move 
and has an orientation induced from $c$, 
then $\mathrm{lk}(K, c') = \mathrm{lk}(K, c) + \varepsilon m$ where $\varepsilon = \pm 1$. 
\item
Assume that $c'$ is
obtained from $c$ by a single $m$--move. 
We give $K$ and $\alpha_m$ parallel orientations, 
$c'$ an orientation induced from $\alpha_m$, 
and $c$ an orientation induced from $c'$.  
Then, an $n$--framing of $c$ becomes
an $(n + 2\mathrm{lk}(K, c) +m)$--framing of $c'$ after an isotopy
in $K(m)$. 
\end{enumerate}
\end{proposition}

We generalize $m$--moves to pairs of seiferters. 

\begin{definition}[$m$--moves to pairs]
\label{m-move for pairs}
Let $c_1 \cup c_2$ be a link in $S^3 - N(K)$. 
Let $\alpha_m \subset \partial N(K)$ be a simple closed curve 
representing a slope $m$, and $b$ a band connecting 
$c_1$ and $\alpha_m$ with $b \cap c_2 =\emptyset$. 
Isotoping the band sum $c_1\,\natural_b\,\alpha_m (\subset S^3 - \mathrm{int}N(K))$ 
away from $\partial N(K)$ without meeting $c_2$, 
we obtain a knot $c'_1 \subset S^3 - N(K)$. 
Then we say that the link $c'_1 \cup c_2$
is obtained from $c_1 \cup c_2$ by 
an \textit{$m$--move} using the band $b$. 
\end{definition}

\begin{proposition}[{\cite[Proposition~2.25]{DMM1}}]
\label{pairband}
Let $K$ be a knot in $S^3$,
and $m$ a slope on $\partial N(K)$. 
Let $c_1 \cup c_2$ and $c'_1 \cup c_2$ be links in 
$S^3 - N(K)$ with each component trivial in $S^3$. 
Suppose that $c'_1 \cup c_2$ is obtained from
$c_1 \cup c_2$ by an $m$--move. 
Then we have: 
\begin{enumerate}
\item
The two ordered links $c_1 \cup c_2$ and $c'_1 \cup c_2$ 
are isotopic in $K(m)$. 
\item 
If $\{c_1, c_2\}$ is a pair of seiferters for $(K, m)$, 
then $\{c'_1, c_2\}$ is also a pair of seiferters for $(K, m)$. 
\end{enumerate}
\end{proposition}

\begin{corollary}[{\cite[Proposition~2.26]{DMM1}}]
\label{m-moves pair}
Let $c_1 \cup c_2$ and $c'_1 \cup c'_2$ be links in 
$S^3 - N(K)$ with each component trivial in $S^3$. 
Let $\alpha_i \subset \partial N(K)$ be a simple closed curve
with slope $m$, 
and $b_i$ a band connecting $c_i$ and $\alpha_i$ such that 
$(c_1 \cup b_1 \cup \alpha_1) \cap (c_2 \cup b_2 \cup \alpha_2)
=\emptyset$. 
Suppose that 
$c'_1 \cup c'_2$ is obtained from 
$(c_1\,\natural_{b_1}\,\alpha_1) \cup (c_2\,\natural_{b_2}\,\alpha_2)$ by 
an isotopy in $S^3 - \mathrm{int}N(K)$.  
Then, $\{ c'_1, c'_2 \}$ is a pair of seiferters for $(K, m)$
if $\{c_1, c_2\}$ is a pair of seiferters for $(K, m)$.  
\end{corollary}

The following proposition will be used to show that
a pair of seiferters for $(K, m)$ does not cobound an annulus 
disjoint from $K$. 

\begin{proposition}[{\cite[Proposition~2.36]{DMM1}}]
\label{irrelevant annular}
Let $c_1$ and $c_2$ be possibly irrelevant seiferters for $(K, m)$ 
with respect to a Seifert fibration 
$\mathcal{F}$ of $K(m)$. 
Suppose that $c_1$ and $c_2$ cobound an annulus $A$ in 
$S^3 - \mathrm{int}N(K)$. 
Then the following hold. 
\begin{enumerate}
\item 
$\mathrm{lk}(c_1, K) = \mathrm{lk}(c_2, K)$. 
\item 
If $K(m)$ is not a lens space, 
then $c_1$ and $c_2$ are regular fibers in $\mathcal{F}$. 
If $K(m)$ is a lens space, 
then we have a Seifert fibration $($possibly distinct
from $\mathcal{F}$$)$ having 
$c_1$ and $c_2$ as regular fibers. 
\end{enumerate}
\end{proposition}

\section{Seiferters for Seifert surgeries on a trefoil knot}
\label{section:seiferter for trefoil}

It is known that there is a hyperbolic knot which 
admits Seifert fibered surgeries
for three successive (integral) surgery slopes.
Well-known examples are the \mbox{$(-1)$--~,} \mbox{$(-2)$--~,} 
$(-3)$--surgeries
on twist knots \cite{BW},
and the \mbox{$17$--~,} \mbox{$18$--~,} $19$--surgeries
on the $(-2, 3, 7)$ pretzel knot \cite{FS}.
In this section, 
we show that for any integer $m \not\in \{-5, -4, -3, -2\}$,
the $m$--~, $(m-1)$--~, $(m-2)$--surgeries on the trefoil knot $T_{-3, 2}$
have a hyperbolic seiferter in common.
Then, arbitrary twists on $T_{-3, 2}$ along the seiferter produce
a knot with three successive Seifert surgeries.
We show that the three successive Seifert fibered surgeries on
twist knots and the $(-2, 3,7)$ pretzel knot arise in this manner.

\begin{figure}[htbp]
\begin{center}
\includegraphics[width=1.0\linewidth]{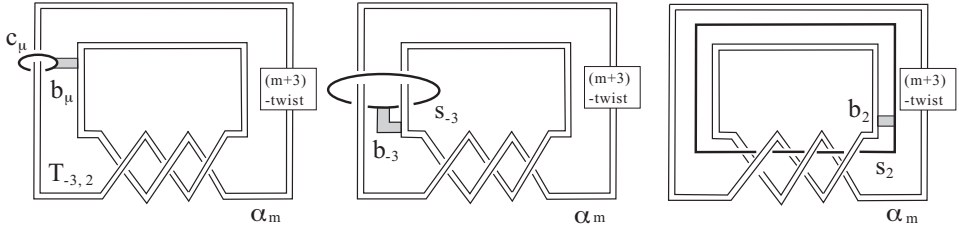}
\caption{Band sums $c_{\mu}\,\natural_{b_{\mu}}\,\alpha_m$,\ 
$s_{-3}\,\natural_{b_{-3}}\,\alpha_m$,\ 
$s_{2}\,\natural_{b_{2}}\,\alpha_m$}
\label{bands}
\end{center}
\end{figure}

Let $\alpha_m \subset \partial N(T_{-3, 2})$ be a simple closed curve
representing a slope $m \in \mathbb{Z}$. 
Let $b_{\mu}$, $b_{-3}$, $b_2$ be the band in
$S^3 -\mathrm{int}N(T_{-3, 2})$ connecting $\alpha$ and
$c_{\mu}$, $s_{-3}$, $s_{2}$, respectively
as described in Figure~\ref{bands}.
We denote by $c_1^m, c_2^m, c_3^m$
the knots obtained from $c_{\mu}, s_{-3}, s_{2}$
by single $m$--moves via these  bands, 
respectively.

\begin{lemma}
\label{seiferter cim}
For any integer $m$ and $i \in \{1, 2, 3\}$, 
$c_i^m$ satisfy the following. 
\begin{enumerate}
\item
The knots $c_1^m$, $c_2^m$, $c_3^m$ are isotopic in $T_{-3, 2}(m)$
to the basic seiferters $c_{\mu}$, $s_{-3}$, $s_2$ for $T_{-3, 2}$,
respectively.
These knots are mutually distinct seiferters for $(T_{-3, 2}, m)$. 
\item
The links $T_{-3, 2} \cup c_1^{m}$, $T_{-3, 2} \cup c_2^{m-1}$, 
and $T_{-3, 2} \cup c_3^{m-2}$ are mutually isotopic in $S^3$.
\end{enumerate}
\end{lemma}

\begin{remark}
\label{trivializing band}
Theorem~1.5 in \cite{IM} shows that
if a band sum of $T_{-3, 2}$ and
its meridian yields a trivial knot,
then such a band is unique up to isotopy.
We thus see that $c_1^m$ is the only seiferter for $(T_{-3 ,2},m)$
obtained from $c_{\mu}$ by an $m$--move.
\end{remark}

\begin{figure}[htbp]
\begin{center}
\includegraphics[width=1.0\linewidth]{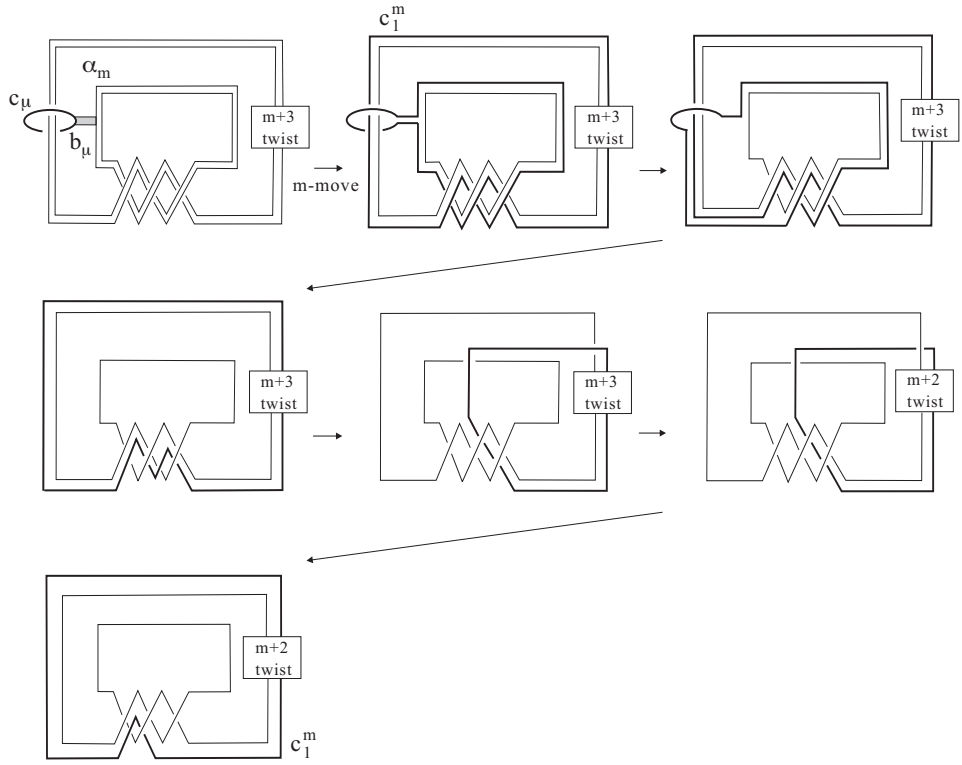}
\caption{$T_{-3, 2} \cup c_1^{m}$}
\label{tc1}
\end{center}
\end{figure}

\begin{figure}[htbp]
\begin{center}
\includegraphics[width=1.0\linewidth]{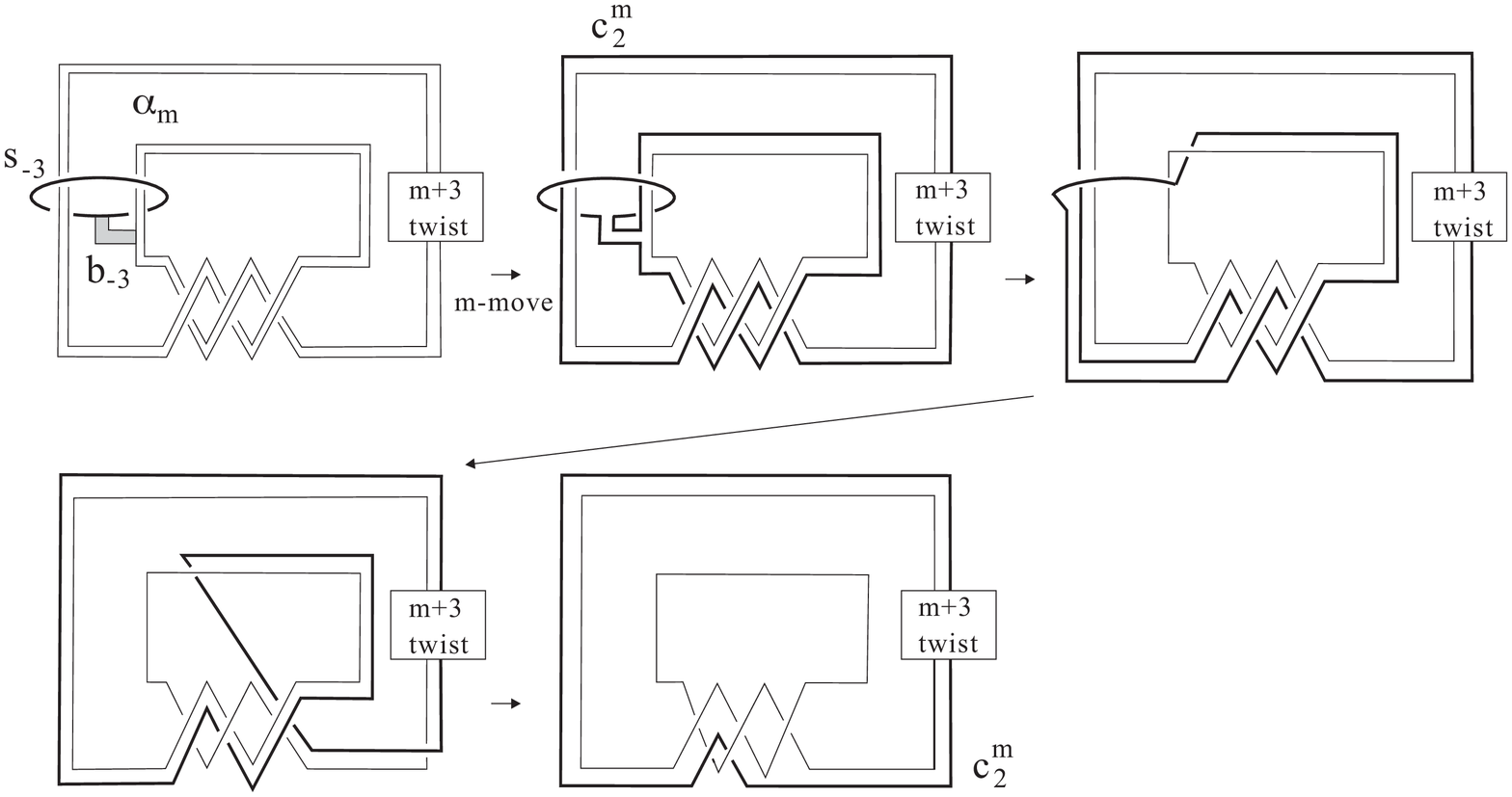}
\caption{$T_{-3, 2} \cup c_2^{m}$}
\label{tc2}
\end{center}
\end{figure}

\begin{figure}[htbp]
\begin{center}
\includegraphics[width=1.0\linewidth]{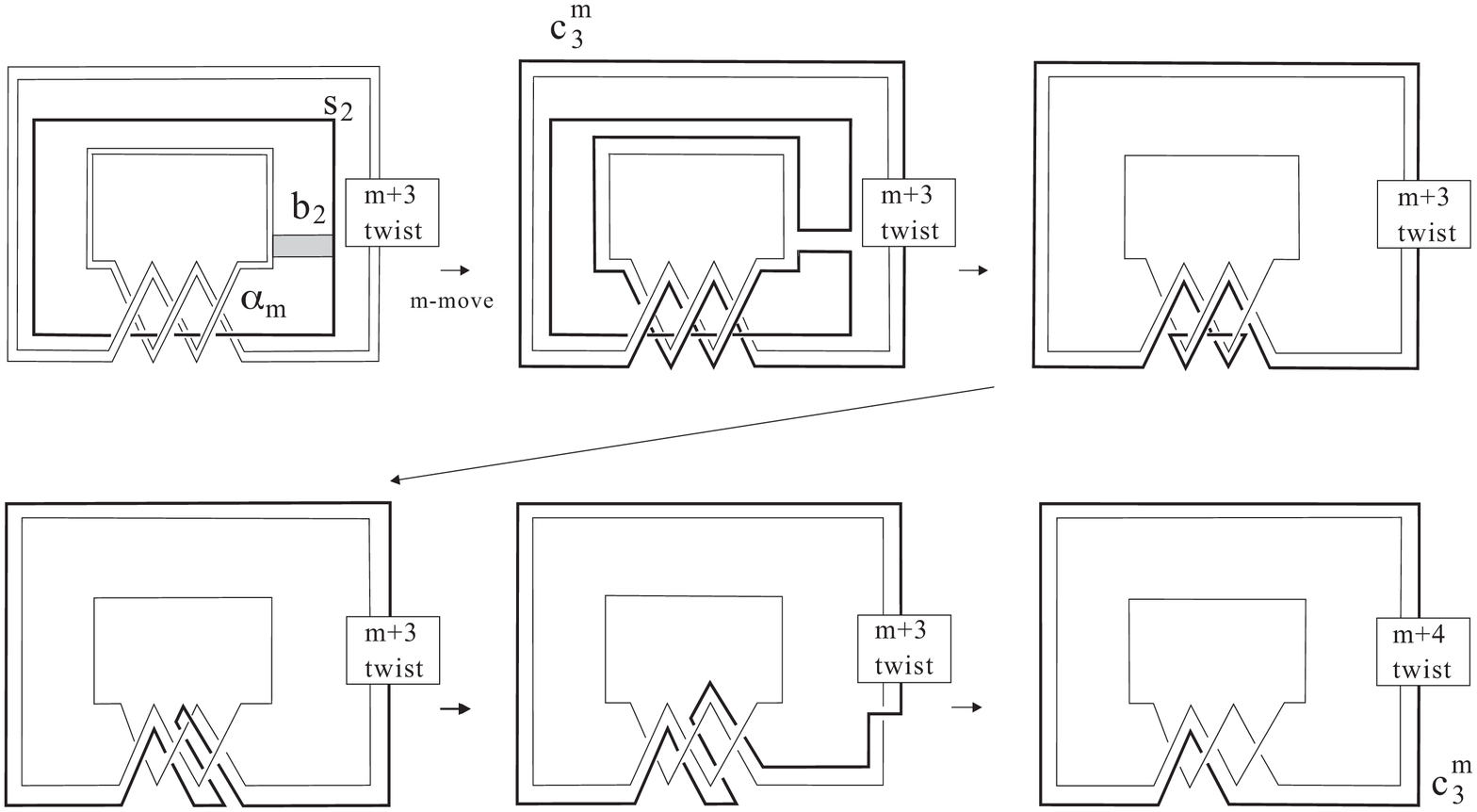}
\caption{$T_{-3, 2} \cup c_3^{m}$}
\label{tc3}
\end{center}
\end{figure}

\noindent
\textit{Proof of Lemma~\ref{seiferter cim}.}
$(1)$ The isotopies in Figures~\ref{tc1}, \ref{tc2}, and \ref{tc3}
show that $c_i^m$ $(i =1, 2, 3)$ are trivial knots in $S^3$. 
Proposition~\ref{move}(1) shows that $c_1^m, c_2^m, c_3^m$
are isotopic to $c_{\mu}, s_{-3}, s_2$, respectively,
and seiferters for $(T_{-3, 2}, m)$. 

We see from Figures~\ref{tc1}, \ref{tc2}, \ref{tc3}
that $|\mathrm{lk} (T_{-3, 2}, c_i^m)| = |m + i|$ $(i=1, 2, 3)$.
These values are mutually distinct except when $m =-2$. 
If $m =-2$, then $|m +i|$ $(i =1, 3)$ are equal.
However, $c_1^{-2}$ and $c_3^{-2}$ are isotopic in $T_{-3, 2}(-2)$
to $c_{\mu}$ and $s_2$, respectively.
Since $c_{\mu}$ and $s_2$ are exceptional fibers of
distinct indices $4$ and $2$ in $T_{-3, 2}(-2)$, 
they are not isotopic in $T_{-3, 2}(-2)$.
It follows that
$c_i^m$ $(i =1, 2, 3)$ are three distinct seiferters
for $(T_{-3, 2}, m)$ for any $m$.
This proves $(1)$. \par

$(2)$ Replace $m$ with $m-1$ (resp. $m-2$) in Figure~\ref{tc2} (resp. Figure~\ref{tc3}), 
we see that $T_{-3,2} \cup c_2^{m-1}$ (resp.\
$T_{-3,2} \cup c_3^{m-2}$) is the same link as $T_{-3,2} \cup c_1^m$.
\QED{Lemma~\ref{seiferter cim}}

Following Lemma~\ref{seiferter cim}(2), 
we denote the seiferters $c_1^{m}, c_2^{m-1}, c_3^{m-2}$ by $c^m$. 
Then Lemma~\ref{seiferter cim} is rephrased as follows. 

\begin{corollary}
\label{seiferter cm}
The knot $c^m$ in $S^3 -T_{-3, 2}$ satisfies the following. 
\begin{enumerate}
\item 
$c^m$ is the seiferter $c^m_1$ for $(T_{-3, 2}, m)$,
$c^{m-1}_2$ for $(T_{-3, 2}, m-1)$, and 
$c^{m-2}_3$ for $(T_{-3, 2}, m-2)$. 
\item
$c^m (= c_1^m),\ c^{m+1} (= c_2^m)$ and 
$c^{m+2} (= c_3^m)$ are mutually distinct seiferters
for $(T_{-3, 2}, m)$. 
\end{enumerate}
\end{corollary}

Since $|\mathrm{lk}(T_{-3, 2}, c^m)| = |m + 1|$,
by twisting $(T_{-3, 2}, m +1 -i)$ $(i = 1, 2, 3)$
$n$ times along $c^m = c^{m +1 -i}_i$
we obtain Proposition~\ref{Knm} below. 
We denote by $K_n^m$
the image of $T_{-3, 2}$ under $n$--twist along $c^m$.
As usual, 
we continue to denote
the image of $c^m$ after twisting along $c^m$ by
the same symbol $c^m$.

\begin{proposition}
\label{Knm}
Let $m$ and $n$ be arbitrary integers.
Then $(K_n^m, m+1 -i +n(m+1)^2)$, 
where $i=1, 2, 3$, 
are Seifert surgeries for which 
$c^m$ is a seiferter.
\end{proposition}

\begin{figure}[htbp]
\begin{center}
\includegraphics[width=0.25\linewidth]{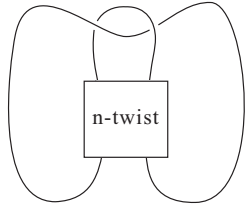}
\caption{Twist knot $Tw(n)$}
\label{twistknot}
\end{center}
\end{figure}

Seifert surgeries given in Proposition~\ref{Knm} contain
three successive Seifert fibered surgeries on
twist knots $Tw(n)$ of Figure~\ref{twistknot},
and the $(-2, 3, 7)$ pretzel knot $P(-2, 3, 7)$.

\begin{proposition}
\label{P-237}
\begin{enumerate}
\item
$(n-1)$--twist along the seiferter $c^{-1}$ converts 
$(T_{-3, 2}, -1)$, $(T_{-3, 2}, -2)$, $(T_{-3, 2}, -3)$
to $(Tw(n), -1), (Tw(n), -2)$, and $(Tw(n), -3)$, respectively. 
\item
$1$--twist along the seiferter $c^{-6}$ converts 
$(T_{-3, 2}, -6)$, $(T_{-3, 2}, -7)$, $(T_{-3, 2}, -8)$
to $(P(-2, 3, 7), 19), (P(-2, 3, 7), 18)$, and $(P(-2, 3, 7), 17)$, respectively. 
\end{enumerate}
\end{proposition}

\begin{figure}[htbp]
\begin{center}
\includegraphics[width=0.9\linewidth]{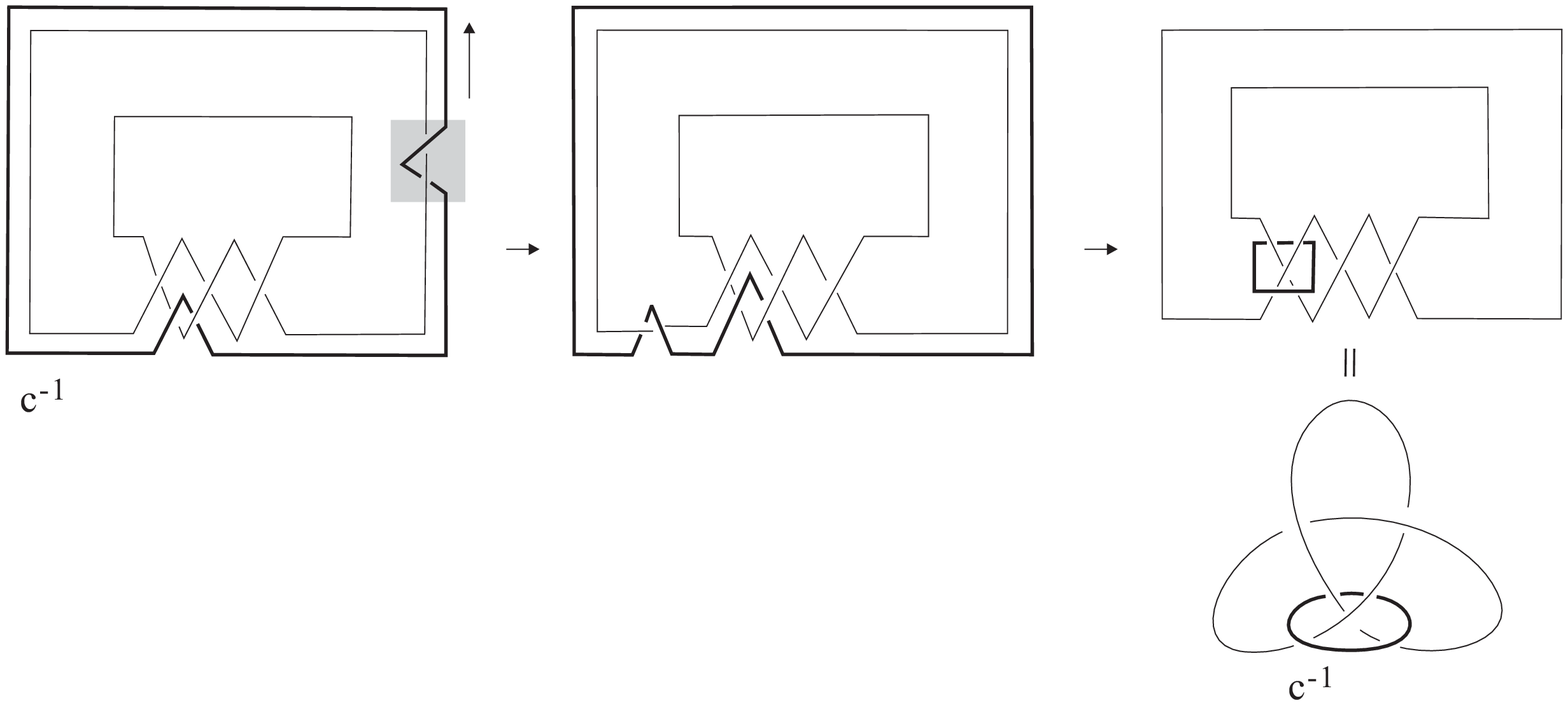}
\vspace{-1em}
\caption{$T_{-3, 2} \cup c^{-1}$}
\label{c-1twist}
\end{center}
\end{figure}

\noindent
\textit{Proof of Proposition~\ref{P-237}.}
$(1)$  We see from Figure~\ref{c-1twist} that
$(n-1)$--twist along $c^{-1}$ converts 
$T_{-3, 2}$ to the twist knot $Tw(n)$,
i.e.\ $K_{n-1}^{-1} =Tw(n)$. 
Since $\mathrm{lk}(c^{-1}, T_{-3, 2}) =0$, 
the surgery slopes do not change under the twisting. 

$(2)$ The sequence of isotopies
in Figures~\ref{P-237c-6_1}, \ref{P-237c-6_2}, and \ref{P-237c-6_3}
shows that $K_1^{-6}$ is the $(-2, 3, 7)$ pretzel knot, 
i.e.\ $1$--twist along $c^{-6}$ converts $T_{-3, 2}$ to $P(-2, 3, 7)$.  
The surgery slopes $17, 18, 19$ are obtained from
Proposition~\ref{Knm} by setting $m = -6,\, n =1$.
\QED{Proposition~\ref{P-237}}

\begin{figure}[htbp]
\begin{center}
\includegraphics[width=0.9\linewidth]{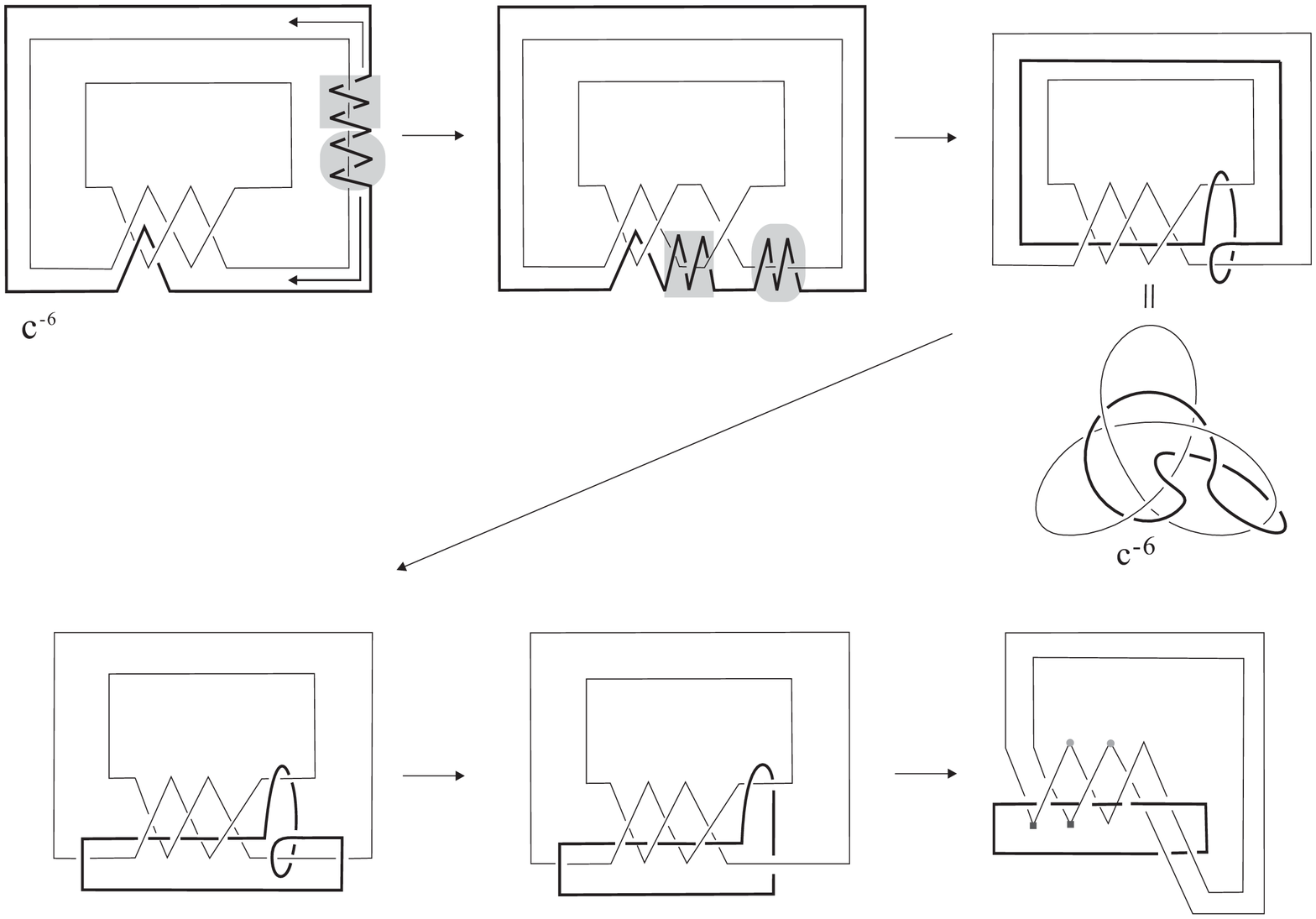}
\caption{Isotopy of $T_{-3, 2} \cup c^{-6}$}
\label{P-237c-6_1}
\end{center}
\end{figure}

\begin{figure}[htbp]
\begin{center}
\includegraphics[width=0.9\linewidth]{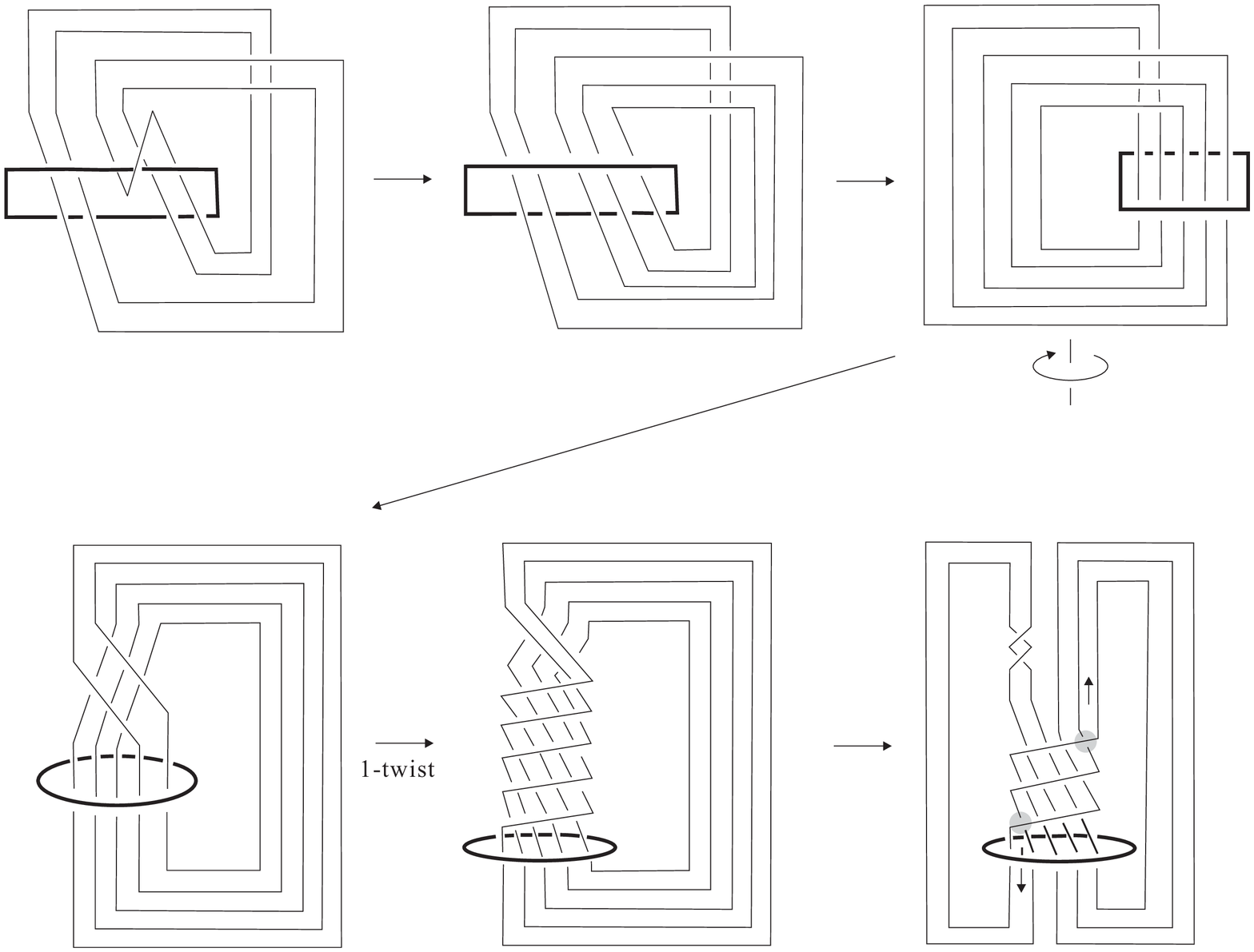}
\caption{Continued from Figure~\ref{P-237c-6_1}}
\label{P-237c-6_2}
\end{center}
\end{figure}

\begin{figure}[htbp]
\begin{center}
\includegraphics[width=0.9\linewidth]{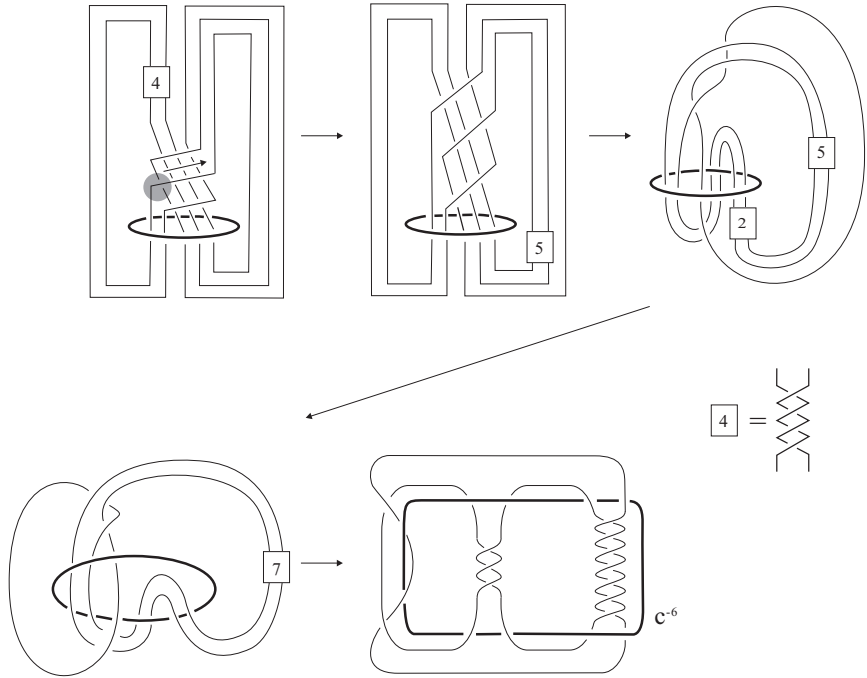}
\caption{Continued from Figure~\ref{P-237c-6_2}:
$1$--twist of $T_{-3,2}$ along $c^{-6}$}
\label{P-237c-6_3}
\end{center}
\end{figure}

\begin{figure}[htbp]
\begin{center}
\includegraphics[width=1.0\linewidth]{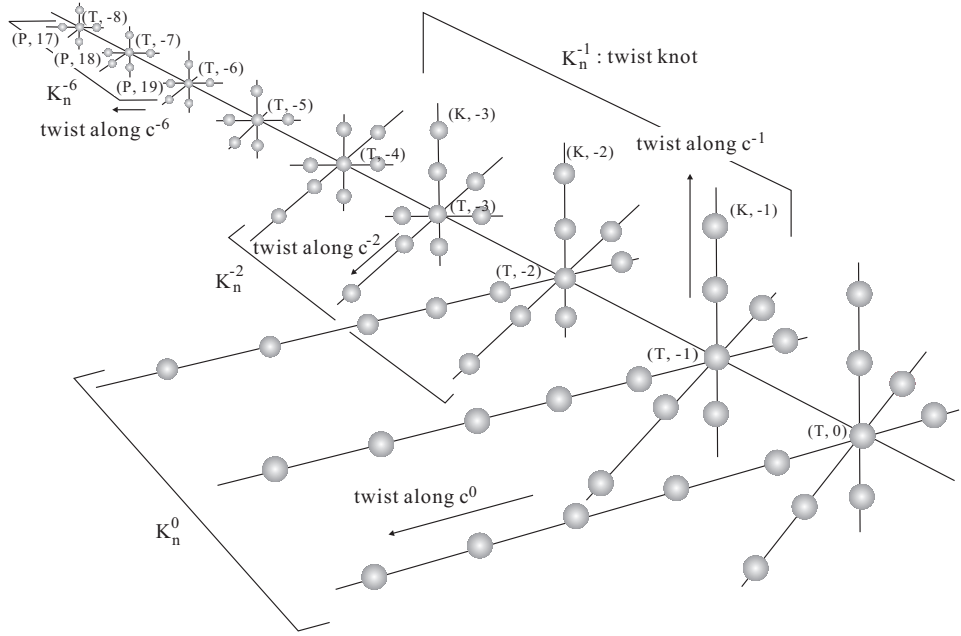}
\caption{Subnetwork generated by
$c^m, c^{m+1}, c^{m+2}$,
where $T = T_{-3, 2}$, $P= P(-2, 3, 7)$, and
$K$ is the figure-eight knot.}
\label{portiontriple}
\end{center}
\end{figure}

Figure~\ref{portiontriple} illustrates
the subnetwork generated by twists along
the seiferters $c^m, c^{m+1}, c^{m+2}$ for $(T_{-3, 2}, m)$.

We calculate the indices of the exceptional fibers of
$K^m_n(m +1-i +n(m +1)^2)$ where $i =1, 2, 3$,
and obtain Proposition~\ref{base orbifold} below.
This proposition will be used to determine when
$K_n^m$ is hyperbolic.

\begin{proposition}
\label{base orbifold}
For any integers $m,n$ the following hold.
\begin{enumerate}
\item
$K^m_n(m +n(m +1)^2)$ is a Seifert fiber space
over the base orbifold
$S^2( 2, 3, |n(m +2)(m +6) +m +n +6| )$.

\item
$K^m_n(m-1 +n(m +1)^2)$ is a Seifert fiber space
over the base orbifold
$S^2( 2, |m +5|, |3n(m +3) -2n +3| )$.

\item
$K^m_n(m-2 +n(m +1)^2)$ is a Seifert fiber space
over the base orbifold
$S^2( 3, |m +4|, |2n(m +4) -3n +2| )$.
\end{enumerate}
\end{proposition}

\noindent
\textit{Proof of Proposition~\ref{base orbifold}.}
The Seifert fiber spaces $K^m_n(m+1 -i +n(m +1)^2)$ $(i =1, 2, 3)$
are obtained from $T_{-3,2}(m+1-i)$ by $n$--twist along $c^m$.   
The seiferter $c^m =c^{m+1 -i}_i$ is isotopic in $T_{-3,2}(m+1 -i)$
to $c_{\mu},\, s_{-3}$, or $s_2$ according as $i=1,\, 2$, or $3$.
Recall that $T_{-3, 2}(m +1 -i)$ has a Seifert fibration
over $S^2(2, 3, |m+7 -i|)$ in which 
$c_{\mu}$, $s_{-3}$, and $s_2$ are fibers of indices $|m+7-i|$,
$3$, and $2$, respectively.

(1) We let $i=1$ in the paragraph above.
Then $K^m_n(m +n(m +1)^2)$ is obtained from
$T_{-3, 2}(m)$ by a surgery along the exceptional fiber
of index $|m +6|$, and thus has a Seifert fibration
over $S^2(2, 3, x)$ for some $x$.
In $T_{-3, 2}(m)$, $c^m_1$ is isotopic to $c_{\mu}$
and further to the core of the filled solid torus $U$.
We set 
$f : H_1(\partial N(c^m_1)) \to H_1(\partial U)
= H_1(\partial N(T_{-3, 2}))$ 
to be the homomorphism induced by this isotopy.
Let $(\mu_c, \lambda_c)$, $(\mu, \lambda)$, and $(\mu', \lambda')$
be preferred meridian--longitude pairs of
$N(c^m_1)$, $N(c_{\mu})$, and $N(T_{-3, 2})$, respectively.
Since $c^m_1$ is obtained from $c_{\mu}$ by an $m$--move,
Proposition~\ref{move}(3) shows that
$0$--framing of $c_{\mu}$ becomes
$(2\mathrm{lk}(T_{-3,2}, c_{\mu}) +m)$--framing of $c^m_1$
after isotopy;
here, $\mathrm{lk}(T_{-3, 2}, c_{\mu}) = 1$ where
$T_{-3, 2}$ are $c_{\mu}$ are oriented so as to satisfy
the assumption in Proposition~\ref{move}(3).
Hence, the isotopy moving
$N(c^m_1)$ to $N(c_{\mu})$ sends $\mu_c, \lambda_c$ to
$\mu, \lambda -(m +2)\mu$ curves on $\partial N(c_{\mu})$.
There is an annulus in $S^3 -\mathrm{int}N(T_{-3, 2} \cup c_{\mu})$
connecting $\lambda (\subset \partial N(c_{\mu}) )$ and
$\mu' (\subset \partial N(T_{-3, 2}))$.
Since $\mu'$ is a longitude of $U$,
the annulus extends to an annulus $A$ connecting
$c_{\mu}$ and the core of $U$.
Isotope $N(c_{\mu})$ to $U$ along $A$.
Then, $\mu$ is sent to
$\lambda' +m\mu'$ curve (a meridian of $U$),
and $\lambda$ is sent to $-\mu'$ curve; 
in fact, $[\lambda' +m\mu'] \cdot [-\mu']
= -[\lambda']\cdot[\mu'] = [\mu']\cdot[\lambda'] =1$. 
Combining these, we obtain
$f([\mu_c]) =[\lambda'] +m[\mu']$ and
$f([\lambda_c]) = -[\mu'] -(m +2)[\lambda' +m\mu']$.
Hence, the image of the $(-\frac{1}{n})$--surgery slope on 
$\partial N(c^m_1)$ is $f([-n \lambda_c + \mu_c])
= n([\mu'] +(m +2) [\lambda' + m\mu']) + [\lambda' +m\mu']
= (n(m +2)+1)[\lambda'] + (mn(m +2) +m +n)[\mu']$.
On the other hand,
a regular fiber on $\partial N(T_{-3, 2})$ represents
$[\lambda'] -6[\mu']$.
Then, the index $x$ is 
$| f([\lambda_c -n\mu_c])\cdot ([\lambda'] -6[\mu']) |$.
Computing this gives the claimed result.

(2) $K^m_n(m -1 +n(m +1)^2)$ is obtained from
$T_{-3, 2}(m -1)$ by a surgery along the exceptional fiber $s_{-3}$
of index $3$, and thus has a Seifert fibration
over $S^2(2, |m +5|, y)$ for some $y$.
In $T_{-3, 2}(m -1)$, $c^{m -1}_2$ is isotopic to $s_{-3}$.
We set $f :H_1(\partial N(c^{m -1}_2)) \to H_1( \partial N(s_{-3}) )$
to be the homomorphism induced by this isotopy.
Let $(\mu_c, \lambda_c)$, $(\mu, \lambda)$
be preferred meridian--longitude pairs of
$N(c^{m -1}_2)$, $N(s_{-3})$, respectively.
Since $c^{m -1}_2$ is obtained from $s_{-3}$ by an $(m -1)$--move,
by Proposition~\ref{move}(3)
$0$--framing of $s_{-3}$ becomes
$(2 \mathrm{lk}(T_{-3, 2}, s_{-3}) +m -1)$--framing of $c^{m-1}_2$
after isotopy;
here, $\mathrm{lk}(T_{-3, 2}, c^{m -1}_{-3}) = 2$.
The reverse isotopy gives $f([\mu_c]) = [\mu]$ and
$f([\lambda_c]) = [\lambda] -(m +3)[\mu]$.
The image of the $(-\frac{1}{n})$--surgery slope on 
$\partial N(c^{m -1}_2)$ is $f([-n\lambda_c +\mu_c])
= -n[\lambda] +( n(m +3) +1 )[\mu]$.
On the other hand,
a regular fiber on $\partial N(s_{-3})$ represents
$-3[\lambda] +2[\mu]$.
Then, the index $y$ is 
$| f([-n\lambda_c + \mu_c])\cdot (-3[\lambda] +2[\mu]) |$.
Assertion~(2) follows from computation.

(3) The proof proceeds in the same manner as in (2)
by replacing $m -1$ with $m -2$,
$c^{m -1}_2$ with $c^{m -2}_3$, and
$s_{-3}$ with $s_2$.
Note that $K^m_n(m -2 +n(m +1)^2)$ is obtained from
$T_{-3, 2}(m -2)$ by a surgery along the exceptional fiber $s_2$
of index $2$, and thus has a Seifert fibration
over $S^2(3, |m +4|, z)$ for some $z$.
In $T_{-3, 2}(m -2)$, $c^{m -2}_3$ is isotopic to $s_2$.
Let $f: H_1(\partial N(c^{m -2}_3)) \to H_1( \partial N(s_2) )$
be the homomorphism induced by the isotopy.
Since $\mathrm{lk}(T_{-3, 2}, s_2) =3$,
Proposition~\ref{move}(3) implies that $f( [\mu_c] ) = [\mu]$ and
$f( [\lambda_c] ) = [\lambda] -(m +4)[\mu]$.
The image of the $(-\frac{1}{n})$--surgery slope on 
$\partial N(c^{m -2}_3)$ is $f([-n\lambda_c + \mu_c])$.
A regular fiber on $\partial N(s_2)$ represents
$2[\lambda] -3[\mu]$.
Then, $z = | f([-n\lambda_c +\mu_c])\cdot ( 2[\lambda] -3[\mu] ) |$.
Computation gives the claimed result.
\QED{Proposition~\ref{base orbifold}}

We determine when $c^m$ is a hyperbolic seiferter
and $K_n^m$ is a hyperbolic knot.

\begin{proposition}
\label{cm is nonbasic}
The link $T_{-3, 2} \cup c^m$ is a hyperbolic link in $S^3$
if and only if $m \not\in \{-5, -4, -3, -2\}$.
In fact, $c^{-4}$, $c^{-3}$, $c^{-2}$ are the same as
the basic seiferters $s_2$, $s_{-3}$, $c_{\mu}$, respectively;
$c^{-5}$ is the $(-1, 2)$ cable of $s_{-3}$.
In \cite[Corollary 3.15(2)]{DMM1},
a $(1, \frac{p +\varepsilon}{2})$ cable of the basic seiferter $s_p$
for $T_{p, 2}$ is called $s_{p, \varepsilon}$,
where $\varepsilon = \pm1$.
Thus $c^{-5}$ is $s_{-3, -1}$.
\end{proposition}

\noindent
\textit{Proof of Proposition~\ref{cm is nonbasic}.}
We first observe that 
$T_{-3, 2} \cup c^m = T_{-3, 2} \cup c_1^m$ in
the last picture of Figure~\ref{tc1} is
isotopic to the Montesinos link
$M(-\frac{1}{2}, \frac{2}{3}, \frac{1}{2m+4})$
given in Figure~\ref{montesinos} below.

\begin{figure}[htbp]
\begin{center}
\includegraphics[width=0.8\linewidth]{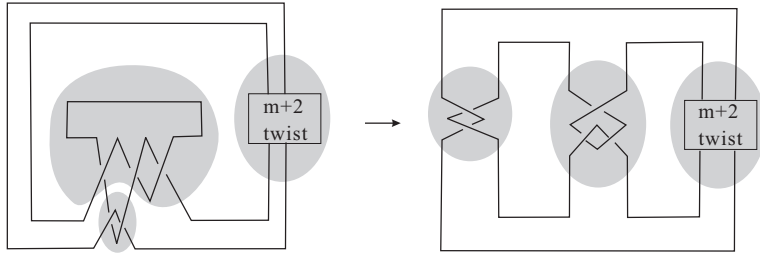}
\caption{$T_{-3, 2} \cup c^m$ is isotopic to
$M(-\frac{1}{2}, \frac{2}{3}, \frac{1}{2m +4})$.}
\label{montesinos}
\end{center}
\end{figure}

If $m = -2$, 
then $c^{-2}$ is a meridian of $T_{-3, 2}$ and 
the exterior $S^3 - \mathrm{int}N( T_{-3, 2} \cup c^{-2} )$ contains
an essential torus. 
If $2m + 4 \ne 0$, 
then the Montesinos link $T_{-3, 2} \cup c^m$ is formed by three rational tangles.  

Assume that $T_{-3, 2} \cup c^m$ is a Seifert link, 
i.e.\ the exterior admits a non-degenerate Seifert fibration. 
Then it turns out that $c^m$ is a non-meridional basic seiferter
for $T_{-3, 2}$.
If $c^m = s_{-3}$, then $|\mathrm{lk}(T_{-2, 3}, c^m)| =|m +1|$
equals to $2$;
if $c^m = s_2$, then $|m +1| =3$.
In the former case, $m = 1, -3$Cand in the latter $m = 2, -4$.
In fact,
$c^{-3}$ is the basic seiferter $s_{-3}$ for $T_{-3, 2}$.
We show that $c^1$ is not a basic seiferter.
Note that
$T_{-3, 2} \cup c^1 = M(-\frac{1}{2}, \frac{2}{3}, \frac{1}{6})$
and 
$T_{-3, 2} \cup c^{-3} = M(-\frac{1}{2}, \frac{2}{3}, -\frac{1}{2})$.
The 2--fold branched covers of $S^3$ along
these links have distinct base orbifolds $S^2(2, 3, 6)$ and
$S^2(2, 3, 2)$.
For a small Seifert fiber space
its Seifert fibrations over $S^2$ are uniquely determined
up to fiber preserving homeomorphism \cite{Hat2, J}.
Thus $T_{-3, 2} \cup c^{-3}$ is not isotopic to $T_{-3, 2} \cup c^{1}$
in $S^3$.
It follows that $c^1$ is not a basic seiferter for $T_{-3, 2}$.
By the same argument we can check that $c^{-4}$ is
the basic seiferter $s_{2}$ for $T_{-3, 2}$,
and $c^2$ is not a basic seiferter for $T_{-3, 2}$.
Hence, $T_{-3, 2} \cup c^m$ is not a Seifert link if $m \not\in \{-4, -3\}$.

In the following, 
assume that $m \not\in \{-2, -3, -4\}$. 
Then, if $T_{-3, 2} \cup c^{m}$
$= M(-\frac{1}{2}, \frac{2}{3}, \frac{1}{2m+4})$
$= M(\frac{1}{2}, -\frac{1}{3}, \frac{1}{2m+4})$ is not hyperbolic, 
\cite[Corollary~5]{Oertel} shows that it is isotopic to 
$M(\frac{1}{2}, -\frac{1}{3}, -\frac{1}{6})$ or its mirror image.
Comparing the indices of the exceptional fibers in
the 2--fold branched cover along $T_{-3, 2} \cup c^m$ with
that along $M(\frac{1}{2}, -\frac{1}{3}, -\frac{1}{6})$,
we see that $2m +4 = \pm 6$ and so $m =1, -5$.
The 2--fold branched cover along
$T_{-3, 2} \cup c^1 = M(\frac{1}{2}, -\frac{1}{3}, \frac{1}{6})$
and that along $M(\frac{1}{2}, -\frac{1}{3}, -\frac{1}{6})$ are not
homeomorphic because
the Euler numbers of these Seifert fibrations are
$\frac{1}{2} -\frac{1}{3}+ \frac{1}{6} = \frac{1}{3}$
and $\frac{1}{2} -\frac{1}{3} -\frac{1}{6} =0$,
i.e.\ distinct up to sign.
This implies that $T_{-3, 2} \cup c^m$ is hyperbolic if $m =1$.
It follows $m = -5$. 
Figure~1.4 in \cite{Oertel} shows that
the exterior of $T_{-3, 2} \cup c^{-5} =
M(\frac{1}{2}, -\frac{1}{3}, -\frac{1}{6})$ in $S^3$
contains an essential torus, 
and furthermore
$c^{-5}$ is the $(-1, 2)$ cable of $s_{-3}$ for $T_{-3, 2}$.
\QED{Proposition~\ref{cm is nonbasic}}

\begin{corollary}
\label{hyperbolic seiferters}
\begin{enumerate}
\item
If $m \le -8$ or $-1\le m$, then $c^m$, $c^{m+1}$, and $c^{m+2}$
are distinct hyperbolic seiferters for $(T_{-3,2}, m)$.
\item
For any integer $m\ne -4$,
$(T_{-3, 2}, m)$ has a hyperbolic seiferter.
\end{enumerate}
\end{corollary}

\noindent
\textit{Proof of Corollary~\ref{hyperbolic seiferters}}.
Proposition~\ref{cm is nonbasic} gives Table~\ref{table of c^m} below.
This table implies assertion~(1).
The table shows that $(T_{-3,2}, m)$ has a hyperbolic seiferter
if $m \not\in \{-4, -5\}$.
For $m = -5$, it is shown in \cite{DMM3} that
the lens surgery $(T_{-3, 2}, -5)$ has a hyperbolic seiferter. 
\QED{Corollary~\ref{hyperbolic seiferters}}

\begin{table}[!h]
\caption{``h" means a hyperbolic seiferter}
\label{table of c^m}
\renewcommand{\arraystretch}{1.2}
\begin{tabular}{|c|c|c|c|c|c|c|c|c|}
\hline 
{\small
\backslashbox{seiferter\\ for $(T_{-3, 2}, m)$}
{\vspace*{-1em}\hspace*{-2em}$m$}
}
   & $m \le -8$& $-7$ & $-6$ & $-5$ & $-4$ & $-3$& $-2$ &$-1 \le m$ \\
\hline
$c^m = c_1^m$  & h & h & h & $s_{-3, -1}$ & $s_2$ & $s_{-3}$ & $c_{\mu}$ & h\\
 \hline
$c^{m+1}= c_2^m$ &h & h & $s_{-3, -1}$ & $s_2$  & $s_{-3}$ & $c_{\mu}$ & h &h\\
 \hline
$c^{m+2}= c_3^m$ & h& $s_{-3, -1}$ & $s_2$& $s_{-3}$ & $c_{\mu}$ & h   &  h &h
\\
\hline
\end{tabular}
\end{table}

We do not know whether $(T_{-3, 2}, -4)$ has a hyperbolic seiferter
or not.
In Section~\ref{section:annular pair for trefoil}, we show
that it has at least six hyperbolic annular pairs of seiferters.

\begin{proposition}
\label{hyperbolic Knm}
The knot $K_{n}^m$ is a hyperbolic knot in $S^3$
if and only if
$m \not\in \{-5, -4, -3, -2\}$, $n\ne 0$, and $(m, n) \ne (-1, -1)$.
\end{proposition}

\noindent
\textit{Proof of Proposition~\ref{hyperbolic Knm}}.
(1) If $m =-5,\, -4,\, -3,$ or $-2$,
then the torus decomposition pieces of
$X =S^3 -\mathrm{N}(T_{-3, 2} \cup c^m)$ are Seifert fiber spaces.
It follows that the exterior of $K_n^m$,
the result of a Dehn filling of $X$, is not hyperbolic
for any $m, n$.
If $n=0$, $K_0^m$ is $T_{-3, 2}$.
If $m =n =-1$, then $K_{-1}^{-1}$ is the twist knot $Tw(0)$,
a trivial knot (see Figure~\ref{twistknot}).
The ``only if" part of
Proposition~\ref{hyperbolic Knm} is proved.

(2) To prove the ``if" part
we start with a proof of Lemma~\ref{satellite} below.

\begin{lemma}
\label{satellite}
No satellite knot has three successive Seifert fibered surgeries.
\end{lemma}

\noindent
\textit{Proof of Lemma~\ref{satellite}}.
This follows from \cite{MM1, MM2}.
By \cite[Theorem~1.2]{MM2}
a satellite knot which is not cabled exactly once admits
at most two (integral) Seifert fibered surgeries and
the slopes are successive.
Assume that $K$ is an $(r, s)$ cable knot and 
$m$--surgery on $K$ is a Seifert fibered surgery.
Then, if $K(m)$ contains an essential torus,
by \cite[Theorem~1.2]{MM1}
 $K$ is the $(2pq \pm1, 2)$ cable of a torus knot $T_{p, q}$
and $m = 4pq$.
If $K(m)$ contains no essential torus, then 
the proof of Theorem~1.4 in \cite{MM1} shows that
$m =rs \pm1$.
Therefore, 
if $K$ is not a $(2pq \pm1, 2)$ cable of $T_{p, q}$, 
it has exactly two Seifert fibered surgeries; 
if $K$ is a $(2pq \pm1, 2)$ cable of $T_{p, q}$, 
it has exactly three Seifert fibered surgeries 
$4pq$, $4pq \pm 1$, $4pq \pm 3$, which are not successive. 
(In the latter case, 
$K(4pq \pm 2)$ is a connected sum of two lens spaces.) 
\QED{Lemma~\ref{satellite}}

\begin{claim}
\label{Kn-6}
$K_n^{-6}$ is a hyperbolic knot for $n \ne 0$.
\end{claim}

\noindent
\textit{Proof of Claim~\ref{Kn-6}}.
If $m = -6$, then
by Proposition~\ref{base orbifold}
$K_n^{-6}(25n -i)$ has a Seifert fibration over
$S^2(2, 3, |n|)$, $S^2(2, 1, |11n -3|)$, or
$S^2(3, 2, |7n -2|)$ according as $i=6,\, 7$, or $8$.
All indices of the exceptional fibers of
these fibrations are nonzero,
so that $(K_n^{-6}, 25n -i)$ $(i =6, 7, 8)$ are
three successive Seifert fibered surgeries for $n \ne 0$.
Hence, $K^{-6}_n$ $(n \ne 0)$ is not a satellite knot
by Lemma~\ref{satellite}.
We also see that $K_n^{-6}$ is a nontrivial knot
because $K_n^{-6}(25n - 8)$ is not a lens space.

For a nontrivial torus knot $T_{p, q}$,
$T_{p, q}(r)$ $(r \in \mathbb{Z})$ is a lens space
if and only if $r =pq \pm 1$;
$T_{p, q}(pq)$ is a connected sum of two lens spaces. 
Now $K_n^{-6}(25n -7)$ is a lens space, 
so that $K_n^{-6}(25n-8)$ or $K_n^{-6}(25n-6)$ is a connected sum of two lens spaces. 
By the assumption $n \ne 0$, this is impossible. 
Hence $K_n^{-6}$ is not a torus knot,
and thus a hyperbolic knot for $n \ne 0$. 
\QED{Claim~\ref{Kn-6}}

We thus assume that
$m \not\in \{-6, -5, -4, -3, -2\}$, $n\ne 0$, and $(m, n) \ne (-1, -1)$.

\begin{claim}
\label{3 small Sfs}
$(K_n^m, m+1 -i +n(m +1)^2)$ $(i = 1, 2, 3)$ are three successive
small Seifert fibered surgeries.
\end{claim}

Then, Lemma~\ref{satellite} shows that
$K_n^m$ is not a satellite knot or a trivial knot.

\smallskip
\noindent
\textit{Proof of Claim~\ref{3 small Sfs}}.
Proposition~\ref{base orbifold}  shows that
the base orbifolds of $K_n^m(m +1-i +n(m +1)^2)$ is
$S^2(2, 3, a)$,
$S^2(2, |m +5|, b)$, or
$S^2(3, |m +4|, c)$
where $a =|n(m +2)(m +6) +m +n +6|$,
$b = |3n(m +3) -2n +3|$,
$c = |2n(m +4) -3n +2|$, according as $i=1$, $2$, or $3$.
Since $|m +5| \ge 2$ and $|m +4|\ge 2$,
it is sufficient to show that
$a, b, c$ are greater than or equal to 2.  

Let $f(x) = nx^2 +(8n +1)x +13 n +6$ $(x \in \mathbb{R})$;
then $|f(m)| = a$.
If $n >0$,
then the axis $y =-\frac{8n +1}{2n}$ of the parabola $y =f(x)$ 
lies between $-5$ and $-4$.
Hence, $f(m) \ge f(-7) =6n -1 \ge 5$, where $m \le -7, m\ge -1$.
If $n <0$, then the fact that $-4 < -\frac{8n +1}{2n} <-3$ implies that
$f(m) \le f(-1) =6n +5 \le -1$;
the last equality holds only if $n=-1$.
Since $(m, n) \ne (-1, -1)$, these results imply
$a =|f(m)| \ge 2$.

Let $g(x) = 3nx +7n +3$ $(x \in \mathbb{R})$; then $b = |g(m)|$
and the zero point $x =-\frac{7n +3}{3n}$ of the linear function $g(x)$
lies between $-4$ and $-1$.   
We see that $|g(m)| \ge |g(-1)| = |4n+3| \ge 1$,
where $m \le -7, m \ge -1$. 
If $|g(m)| = 1$, then $(m, n) = (-1, -1)$, a contradiction. 
Thus $b \ge 2$. 

Regarding $c = h(x)$ where $h(x) = 2nx +5n+2$,
the zero point $-\frac{5n +2}{2n}$ of $h$ lies between
$-4$ and $-1$.
Thus, for $m \le -7$ and $m\ge -1$ we have
$|h(m)| \ge |h(-1)| = |3n +2| \ge 1$. 
The equality $|h(m)| = 1$ holds only if $(m, n) = (-1, -1)$,
an excluded case.  
Hence, $c \ge2$ as desired.
\QED{Claim~\ref{3 small Sfs}}

Assume for a contradiction
that $K_n^m$ is a torus knot $T_{p, q}$ where $|p| >q \ge2$
for some $m, n$ satisfying our assumption.
For simplicity, set $d =m +1 +n(m +1)^2$.
Then $K_n^m(d -i)$ admits a Seifert fibration
over $S^2(|p|, q, |pq -d +i|)$
for $i \in \{ 1, 2, 3 \}$.
Thus the unordered triples of the indices of exceptional fibers satisfy
$(2, 3, a ) = (|p|, q, |pq -d+1|)$, 
$(2, |m + 5|, b ) = (|p|, q, |pq -d+2|)$, and
$(3, |m + 4|, c ) = (|p|, q, |pq -d +3|)$ for some integers $m, n$,
$|p| >q \ge2$.
Since $|pq - d + i| \ne 0, 1$ for $i = 1, 2, 3$ (Claim~\ref{3 small Sfs}), 
the indices $|pq -d+1|, |pq -d+2|, |pq -d+3|$ are
mutually distinct.
Hence, all triples have 2 and 3 in common,
so that $|m +4| = 2$ or $c =2$.
The former case implies $m =-2, -6$,
a contradiction. 
The latter case $c =|2nm +5n +2| = 2$ implies
$(m, n) =(-3, 4), (-2, -4)$, a contradiction.
Therefore, $K_n^m$ is neither a satellite knot nor
a torus knot,
so that $K_n^m$ is a hyperbolic knot.
\QED{Proposition~\ref{hyperbolic Knm}}

Proposition~\ref{hyperbolic Knm} and its proof imply
the following theorem,
which generalizes a previous result in \cite{MSong}.

\begin{theorem}
\label{triple}
The following (1) and (2) hold.
\begin{enumerate}
\item
For any integer $m$, there is a hyperbolic knot $K$
such that $(K, m)$, $(K, m+1)$, $(K, m+2)$ are
small Seifert fibered surgeries.
\item
If $m \ne -2$, the hyperbolic knot $K$ in (1) above can be chosen
so that three successive surgeries in (1) arise
from three successive Seifert surgeries
on $T_{-3, 2}$ by twisting along a common seiferter.
\end{enumerate}
\end{theorem}

\noindent
\textit{Proof of Theorem~\ref{triple}.}
Claim~\ref{3 small Sfs} shows that
$(K_n^0, n)$, $(K_n^0, n-1)$, $(K_n^0, n-2)$
are small Seifert fibered surgeries, where $n \ne 0$.
By Proposition~\ref{hyperbolic Knm} $K_n^0$ $(n \ne 0)$ is hyperbolic.
Regarding $0$--, $(-1)$--, $(-2)$--surgeries,
take the mirror images of $(K_2^0, 2)$, $(K_2^0, 1)$, $(K_2^0, 0)$.
Then we obtain
$(-2)$--, $(-1)$--, $0$--surgeries on the mirror image $K'$ of $K_2^0$,
small Seifert fibered surgeries on a hyperbolic knot.
Note that these surgeries on $K'$
also arise from three successive surgeries on $T_{-3,2}$ 
after twisting along a common seiferter.
\QED{Theorem~\ref{triple}}

\begin{question}
\label{0,-1,-2}
Is the assumption $m \ne -2$ in Theorem~\ref{triple}(2)
necessary?
\end{question}

\section{Annular pairs for Seifert surgeries on a trefoil knot}
\label{section:annular pair for trefoil}

As observed in Section~\ref{section:seiferter for trefoil}, 
$(T_{-3, 2}, m)$ has six seiferters 
$c_{\mu}, s_{-3}, s_2, c_1^m, c_2^m$, $c_3^m$.
We completely determined which are basic seiferters
or hyperbolic seiferters (Proposition~\ref{cm is nonbasic}). 
In this section, 
we obtain annular pairs of seiferters by placing
any two of these six seiferters in adequate positions. 
We first prove Lemma~\ref{annularT-32} below
by applying $m$--moves to basic annular pairs for $(T_{-3, 2}, m)$.

\begin{figure}[htbp]
\begin{center}
\includegraphics[width=0.9\linewidth]{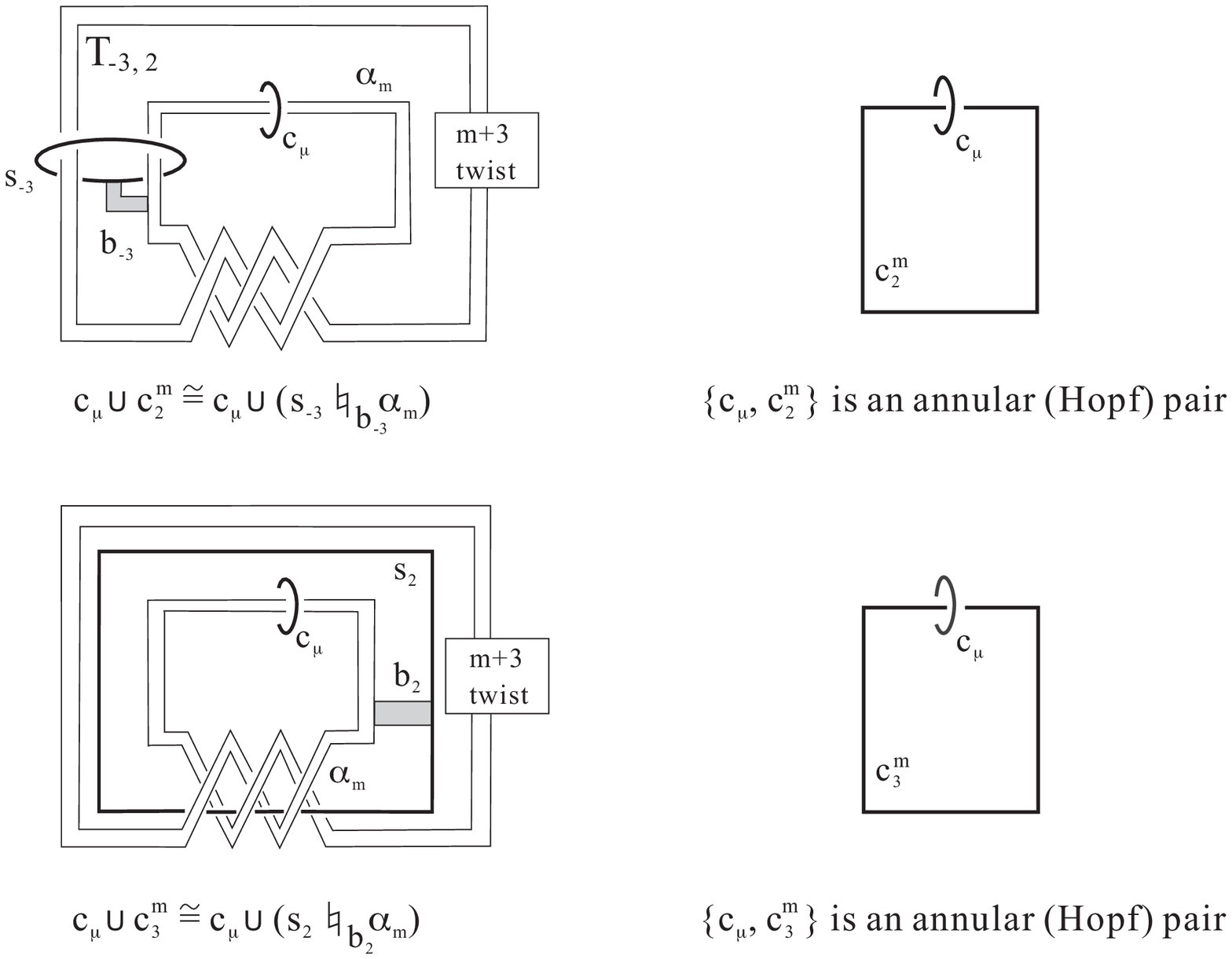}
\caption{Annular pairs of seiferters $\{c_{\mu}, c_2^m \}$,
$\{c_{\mu}, c_3^m \}$ for $(T_{-3, 2}, m)$}
\label{T-32annularpair1}
\end{center}
\end{figure}

\begin{figure}[htbp]
\begin{center}
\includegraphics[width=0.9\linewidth]{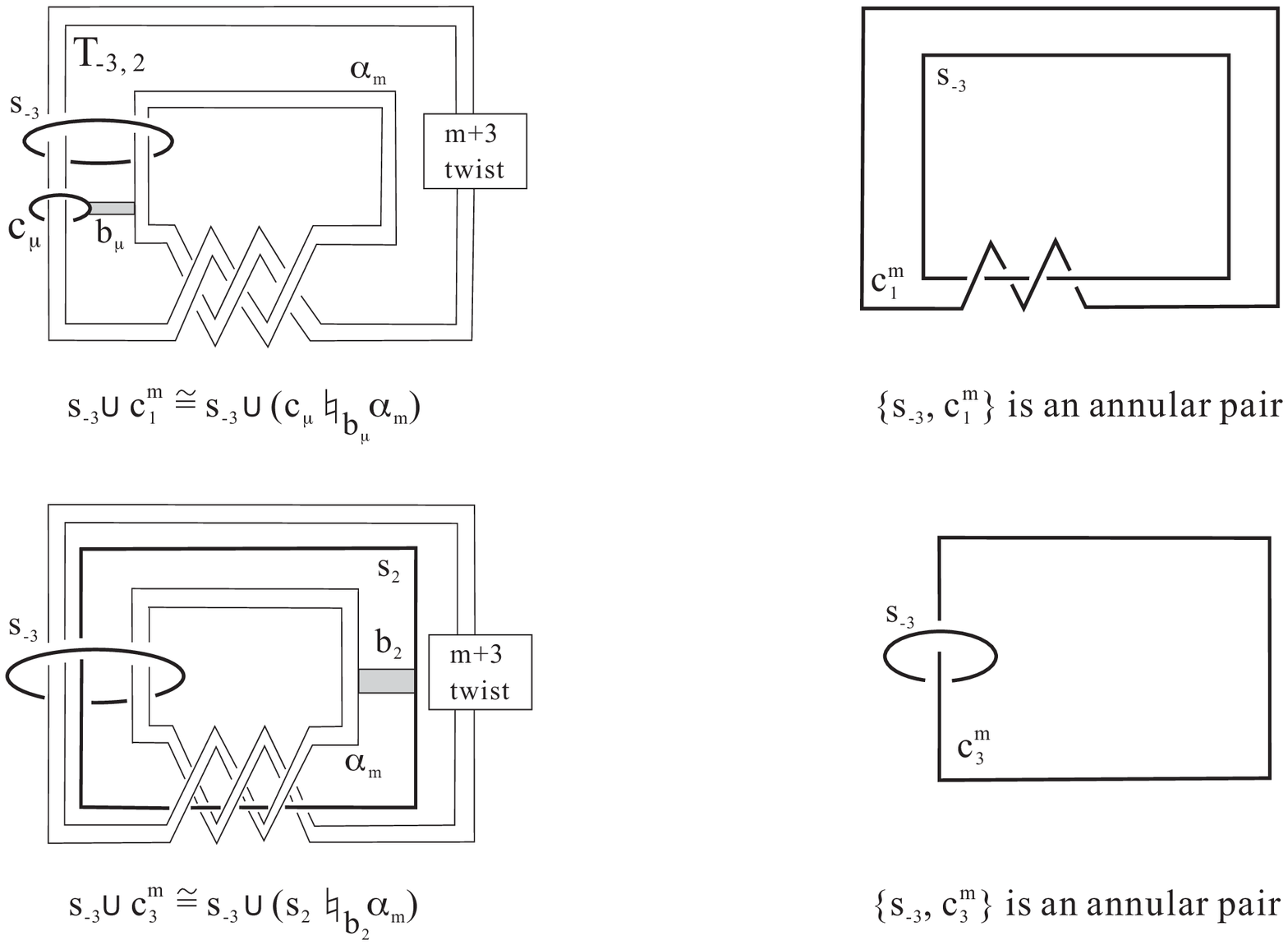}
\caption{Annular pairs of seiferters
$\{ s_{-3}, c_1^m \}$, $\{s_{-3}, c_3^m \}$ for $(T_{-3, 2}, m)$}
\label{T-32annularpair2}
\end{center}
\end{figure}

\begin{figure}[htbp]
\begin{center}
\includegraphics[width=0.9\linewidth]{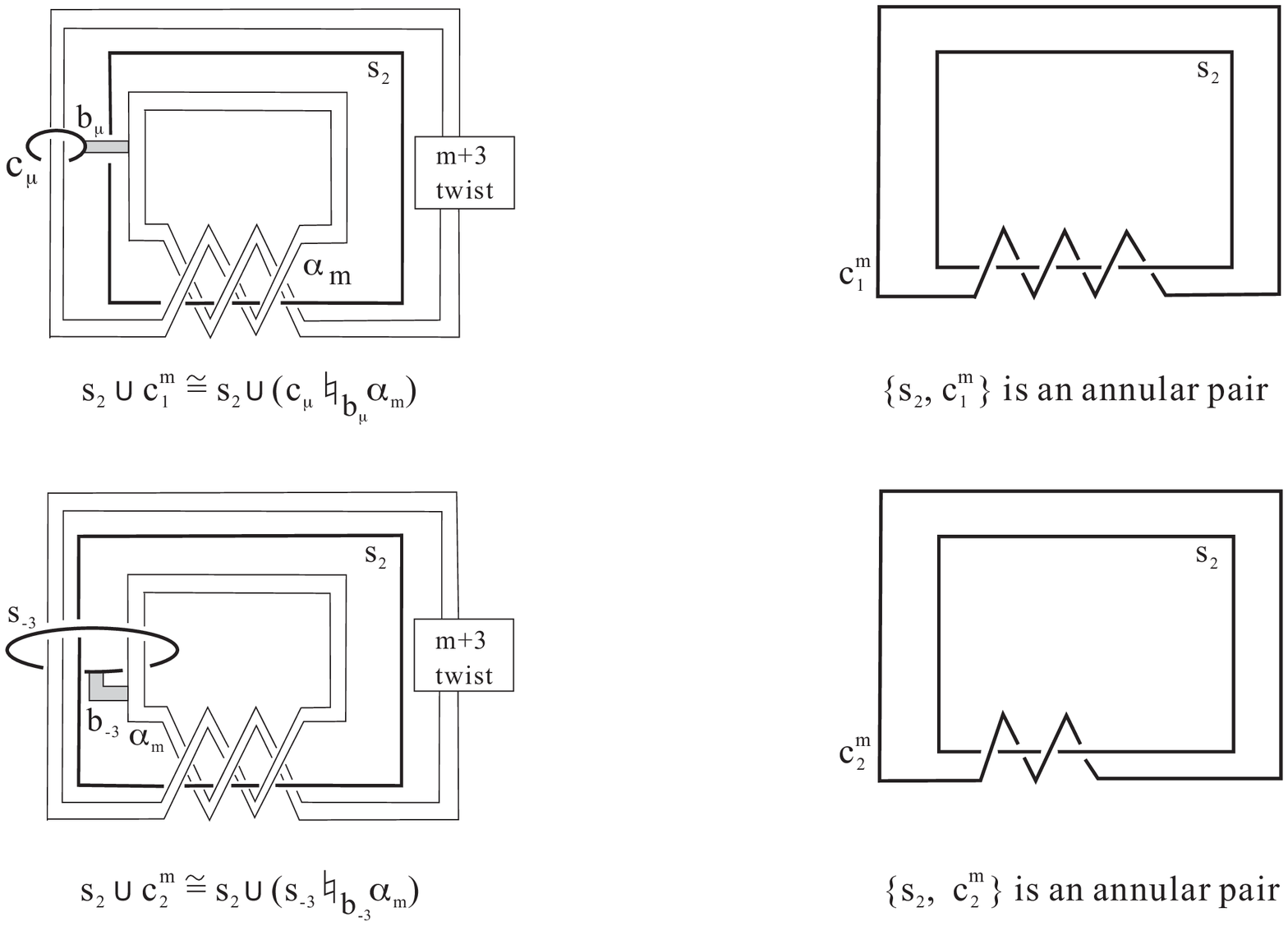}
\caption{Annular pairs of seiferters $\{s_2, c_1^m \}$,
$\{s_2, c_2^m \}$ for $(T_{-3, 2}, m)$}
\label{T-32annularpair3}
\end{center}
\end{figure}

\begin{figure}[htbp]
\begin{center}
\includegraphics[width=0.9\linewidth]{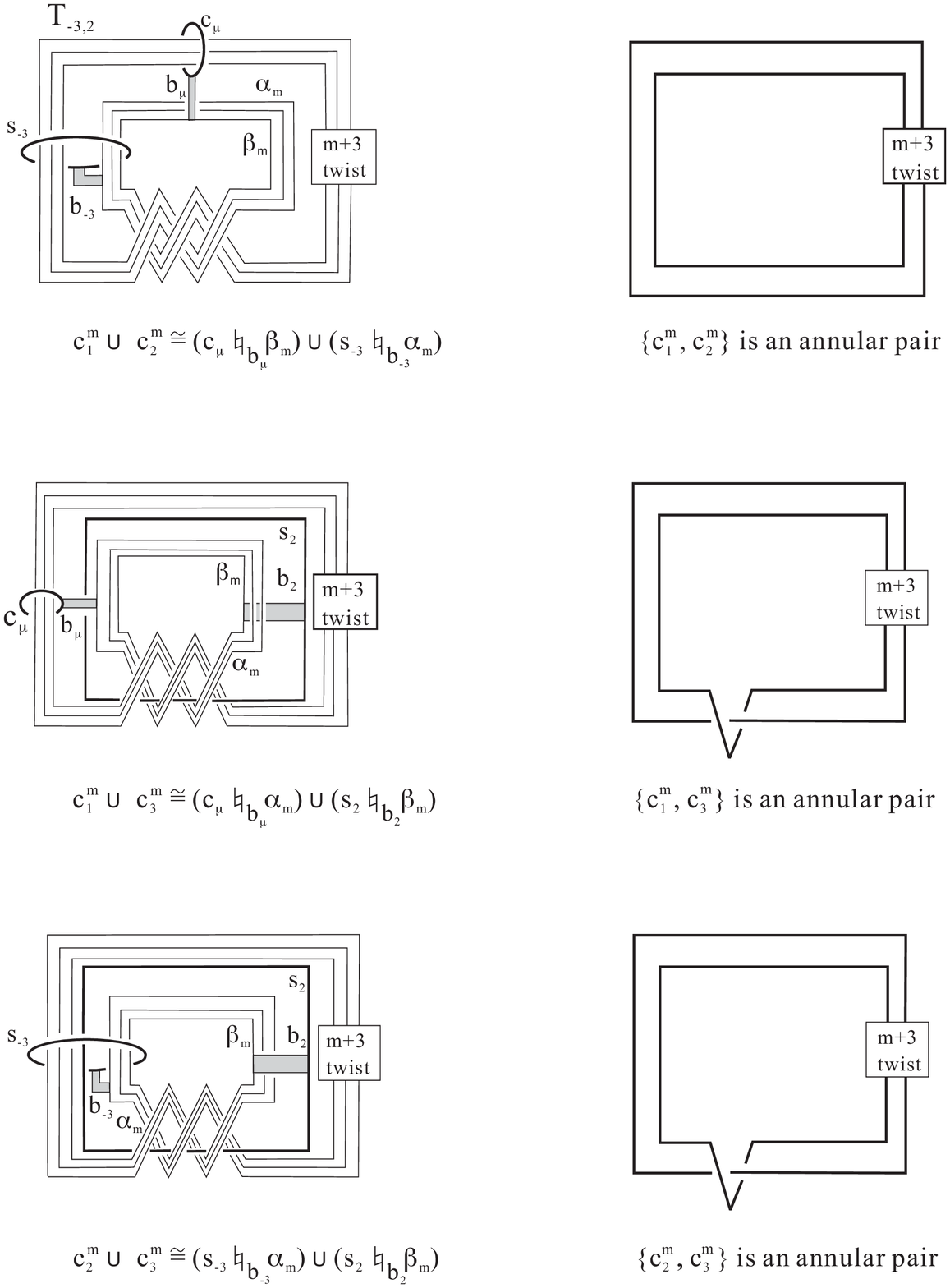}
\caption{Annular pairs of seiferters $\{ c_1^m, c_2^m \}$, 
$\{c_1^m, c_3^m\}$, 
$\{c_2^m, c_3^m\}$ for $(T_{-3, 2}, m)$}
\label{T-32annularpair4}
\end{center}
\end{figure}

\begin{lemma}
\label{annularT-32}
$\{c_{\mu}, c_2^m\}$, $\{c_{\mu}, c_3^m\}$
in Figure~\ref{T-32annularpair1}, 
$\{s_{-3}, c_1^m\}$, $\{s_{-3}, c_3^m\} (\ne \{s_{-3}, c_3^{-5}\})$ 
in Figure~\ref{T-32annularpair2},
$\{s_2, c_1^m\}$, $\{s_2, c_2^m\} (\ne \{s_2, c_2^{-5}\})$ 
in Figure~\ref{T-32annularpair3}, 
$\{c_1^m, c_2^m\}$, $\{c_1^m, c_3^m\}$ and $\{c_2^m, c_3^m\}$ 
in Figure~\ref{T-32annularpair4} 
are annular pairs of seiferters 
for $(T_{-3, 2}, m)$ for any $m$. 
Each of these pairs is isotopic in $T_{-3,2}(m)$
to a basic annular pair for $(T_{-3, 2}, m)$.
\end{lemma}

\begin{remark}
\label{irrelevant cases}
The excluded pairs $\{s_{-3}, c_3^{-5}\}$ and $\{s_2, c_2^{-5}\}$
in Figures~\ref{T-32annularpair3}, \ref{T-32annularpair4} with $m = -5$
cobound annuli in $S^3 - \mathrm{int}N(T_{-3, 2})$. 
By Remark~\ref{rem:irrelevant}  
they are not regarded as annular pairs of seiferters
for $(T_{-3, 2}, -5)$. 
\end{remark}

\noindent
\textit{Proof of Lemma~\ref{annularT-32}.} 
In Figures~\ref{T-32annularpair1}--\ref{T-32annularpair3},
each pair uses one band,
and is obtained from a basic annular pair of seiferters
by an $m$--move.
Then Proposition~\ref{pairband}(2) shows that such a pair is
a pair of seiferters.
Regarding pairs in Figure~\ref{T-32annularpair4},
each of them uses two mutually disjoint bands
and satisfies the assumption in Corollary~\ref{m-moves pair}.
It then follows from Corollary~\ref{m-moves pair}
that each pair in Figure~\ref{T-32annularpair4}
is a pair of seiferters.
\par
The fact that all the pairs cobound annuli in $S^3$ is shown in
Figures~\ref{T-32annularpair1}--\ref{T-32annularpair4}. 
It remains to show that all the pairs 
in Lemma~\ref{annularT-32} do not cobound annuli 
in $S^3 - \mathrm{int}N(T_{-3, 2})$. 
Recall that 
$|\mathrm{lk}(T_{-3, 2}, c_{\mu})| = 1, 
|\mathrm{lk}(T_{-3, 2}, s_{-3})| = 2, 
|\mathrm{lk}(T_{-3, 2}, s_2)| = 3$, and 
$|\mathrm{lk}(T_{-3, 2}, c_i^m)| = |m + i|$. 
Hence, 
if a pair of seiferters
in Figures~\ref{T-32annularpair1}--\ref{T-32annularpair3}
satisfies the condition that its components
have the same linking numbers
with $T_{-3, 2}$ (up to sign), 
then the pair is one of the following list: 
$\{c_{\mu}, c_2^m \}$ with $m = -3, -1$, 
$\{c_{\mu}, c_3^m \}$ with $m = -4, -2$, 
$\{s_{-3}, c_1^m \}$ with $m = -3, 1$, 
$\{s_{-3}, c_3^m\}$ with $m = -5, -1$, 
$\{s_2, c_1^m \}$ with $m = -4, 2$, 
$\{s_2, c_2^m \}$ with $m = -5, 1$. 
Then, Proposition~\ref{irrelevant annular}(1) guarantees that 
the pairs not on the above list
are (relevant) annular pairs of seiferters 
for $(T_{-3, 2}, m)$. 
Proposition~\ref{irrelevant annular}(2) shows that
if $T_{-3, 2}(m)$ is not a lens space and
a pair of seiferters $\{ c_1, c_2 \}$ for $(T_{-3,2},m)$ cobounds
an annulus is $S^3 -T_{-3, 2}$, 
then $c_1$ and $c_2$ are regular fibers in $T_{-3, 2}(m)$.
Each seiferter in a pair listed above is an exceptional fiber
in $T_{-3, 2}(m)$ if $T_{-3, 2}(m)$ is not a lens space.
This is because each basic seiferter is
a (possibly degenerate) exceptional fiber in $T_{-3, 2}(m)$ 
if $T_{-3, 2}(m)$ is not a lens space. 
Hence, 
Proposition~\ref{irrelevant annular}(2) narrows the list above to
the cases when $T_{-3, 2}(m)$ is a lens space: 
$\{ s_{-3}, c_3^{-5} \}$ and $\{ s_2, c_2^{-5} \}$;
both are pairs of seiferters for $(T_{-3, 2}, -5)$. 
\QED{Lemma~\ref{annularT-32}}

\begin{figure}[htbp]
\begin{center}
\includegraphics[width=0.9\linewidth]{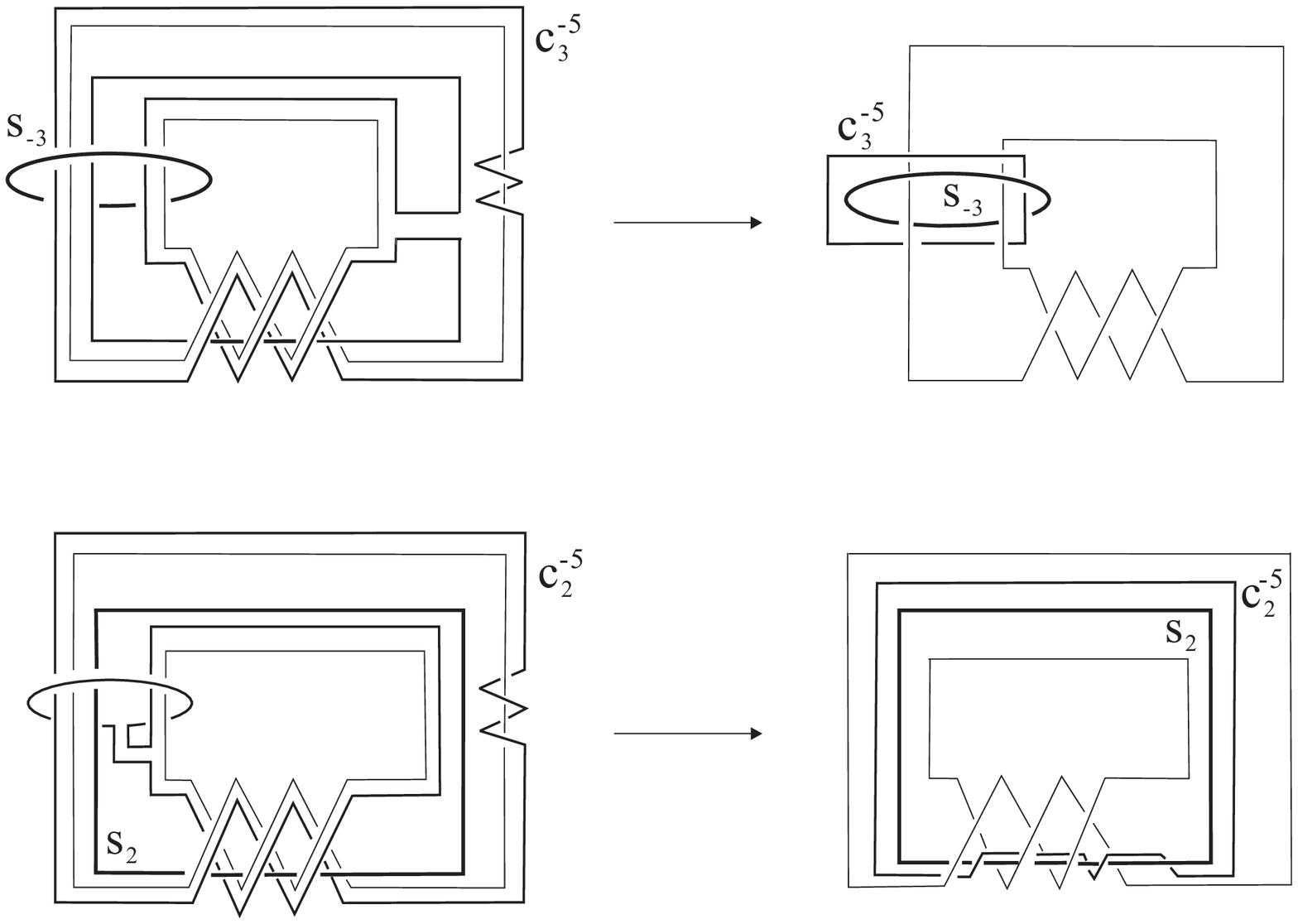}
\caption{$\{s_{-3}, c_3^{-5}\}$ and $\{s_2, c_2^{-5}\}$ are irrelevant annular pairs of seiferters for $(T_{-3, 2}, -5)$.}
\label{coboundannuli}
\end{center}
\end{figure}

Assume $|m +6| >3$.
Then, no matter where the two seiferters $c_{\mu}$ and $c_1^m$
are placed in $S^3 -T_{-3, 2}$,
$\{c_{\mu}, c_1^m\}$ cannot be a pair of seiferters
for $(T_{-3, 2}, m)$ for the following reason.
First note that
any Seifert fibration on $T_{-3, 2}(m)$
has three exceptional fibers with mutually distinct indices
$2, 3, |m +6|$.
Since $c_{\mu}$ and $c_1^m$ are isotopic in $T_{-3, 2}(m)$,
their exteriors in $T_{-3, 2}(m)$ are homeomorphic.
Therefore, if $T_{-2, 3}(m)$ has a Seifert fibration
in which $c_{\mu}$ and $c_1^m$ are fibers simultaneously, 
they are regular fibers.
Then $T_{-2, 3}(m) -\mathrm{int}N(c_{\mu})$ is a Seifert fiber
space over $D^2(2, 3, |m+6|)$.
This contradicts the fact that
$T_{-2, 3}(m) -\mathrm{int}N(c_{\mu})$ is a Seifert fiber space
over $D^2(2, 3)$.
For the same reason, $\{ s_{-3}, c_2^m \}$ and $\{s_2, c_3^m\}$
cannot be pairs of seiferters for $(T_{-3, 2}, m)$, 
where $|m +6| > 3$.
On the other hand,
for $m = -5, -7$, using the flexibility of Seifert fibrations 
on the lens space $T_{-3, 2}(m)$, 
we obtain the following.

\begin{figure}[htbp]
\begin{center}
\includegraphics[width=0.9\linewidth]{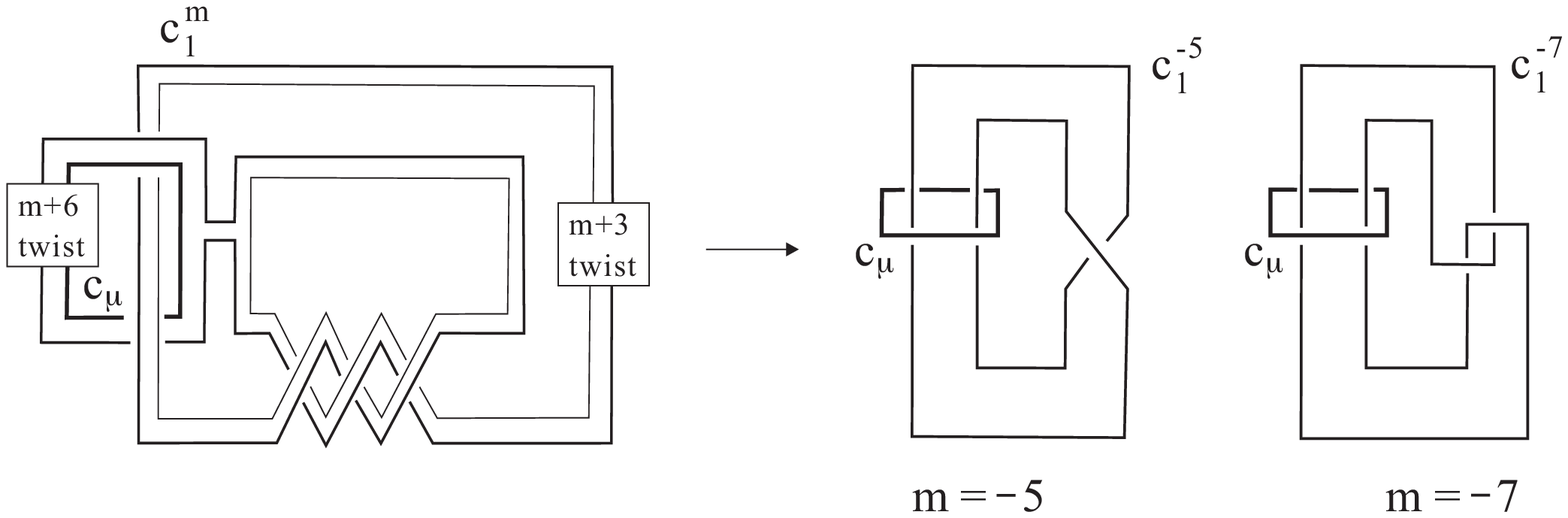}
\caption{Pairs of seiferters $\{c_{\mu}, c_1^m\}$ for $(T_{-3, 2}, m)$,
where $m \in \{-5, -7\}$.}
\label{T-32m-5-7}
\end{center}
\end{figure}

\begin{figure}[htbp]
\begin{center}
\includegraphics[width=1.0\linewidth]{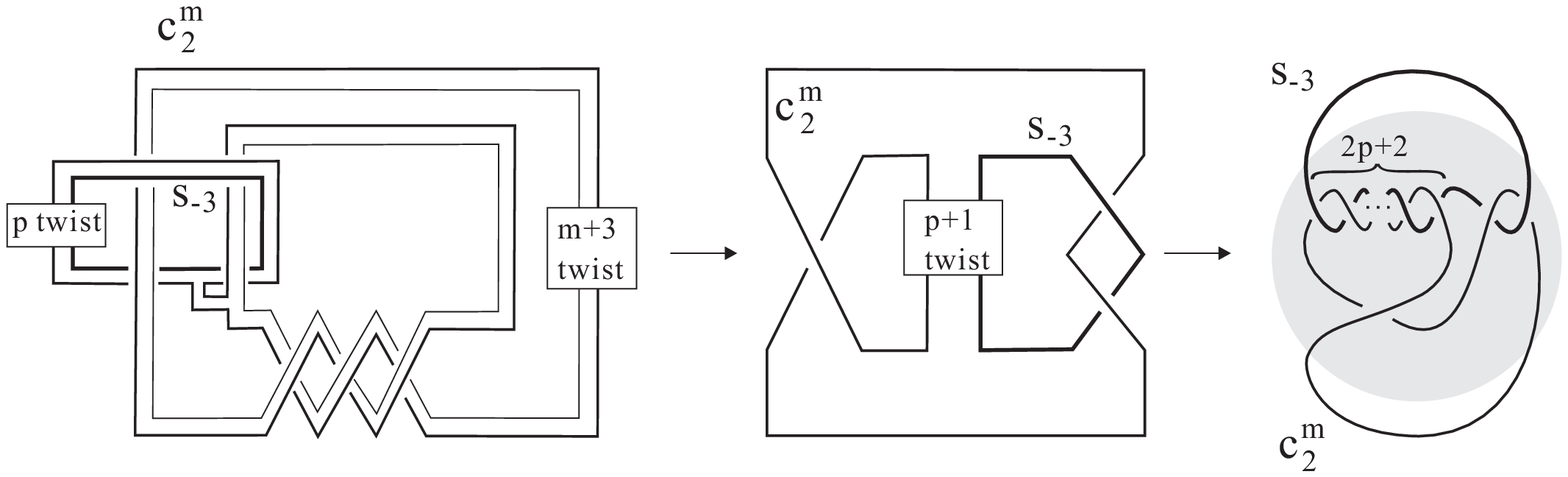}
\caption{Pairs of seiferters $\{s_{-3}, c_2^m\}$ for $(T_{-3, 2}, m)$,
where $m \in \{-5, -7\}$.}
\label{T-32m-5-7_2}
\end{center}
\end{figure}

\begin{figure}[htbp]
\begin{center}
\includegraphics[width=1.0\linewidth]{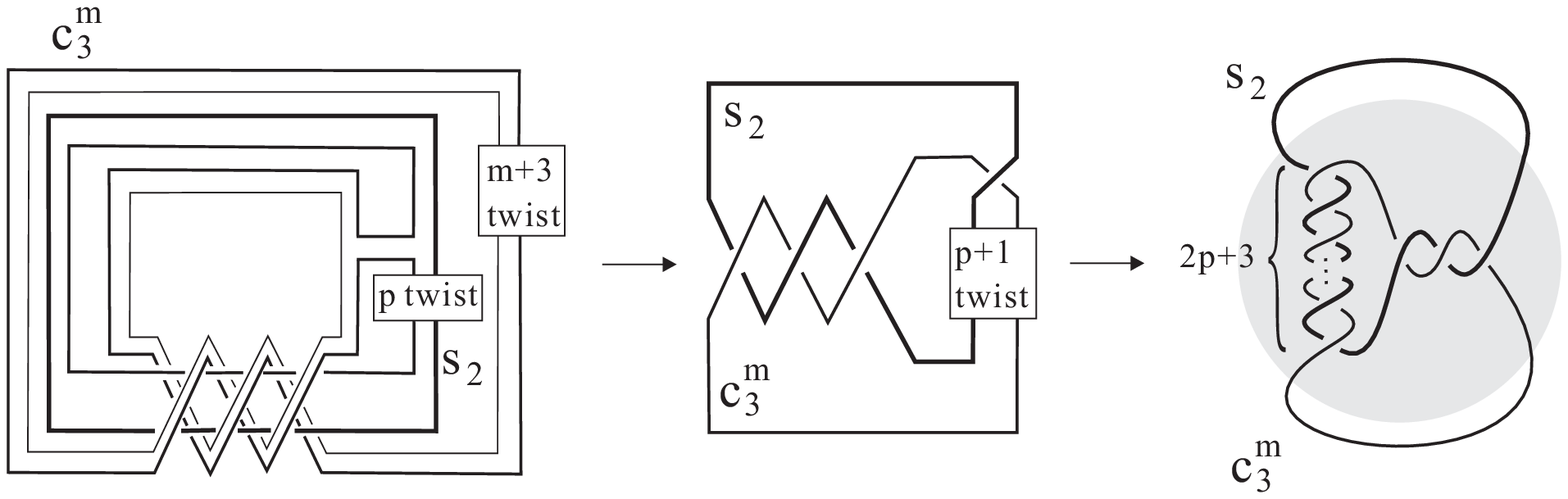}
\caption{Pairs of seiferters $\{s_2, c_3^m\}$ for $(T_{-3, 2}, m)$,
where $m \in \{-5, -7\}$.}
\label{T-32m-5-7_3}
\end{center}
\end{figure}

\begin{lemma}
\label{seifert pair example}
For $m= -5, -7$, assertions~(1), (2), (3) below hold.
Each $($annular$)$ pair of seiferters obtained below consists of fibers
in a non-degenerate Seifert fibration of $T_{-3, 2}(m)$.
\begin{enumerate}
\item
$\{c_{\mu}, c_1^{-5}\}$ in Figure~\ref{T-32m-5-7} is an annular pair of seiferters for 
$(T_{-3, 2}, -5)$; 
$\{c_{\mu}, c_1^{-7}\}$ in Figure~\ref{T-32m-5-7} is a pair of seiferters, 
but not an annular pair of seiferters for $(T_{-3, 2}, -7)$.  
\item
$\{ s_{-3}, c_2^m \}$ in Figure~\ref{T-32m-5-7_2} is 
a pair of seiferters for $(T_{-3, 2}, m)$ for any integer $p$; 
$\{ s_{-3}, c_2^m \}$ is an annular pair of seiferters for $(T_{-3, 2}, m)$
exactly when $p = 0, -1$. 
\item
$\{s_2, c_3^m\}$
in Figure~\ref{T-32m-5-7_3}
is a pair of seiferters for $(T_{-3, 2}, m)$ for any integer $p$; 
$\{s_2, c_3^m\}$ is an annular pair of seiferters for $(T_{-3, 2}, m)$
exactly when $p =-1, -2$.
\end{enumerate} 
\end{lemma}

\noindent
\textit{Proof of Lemma~\ref{seifert pair example}.}
Take $m \in \{-5, -7\}$.

(1) The link $c_{\mu} \cup c_1^m$ in Figure~\ref{T-32m-5-7}
is isotopic in $T_{-3, 2}(m)$ to the union of $c_{\mu}$ and
the $(m+6, 1)$ cable of $N(c_{\mu})$.
Note that since $m+6 = \pm 1$, the $(m+6, 1)$ cable
is the $(1, m+6)$ cable.
For any integer $n$, 
an isotopy in $T_{-3, 2}(n)$ sending $c_{\mu}$ to the core of the filled solid torus
sends the $(1, n +6)$ cable of $c_{\mu}$ to the $(-6, 1)$ cable
of $T_{-3, 2}$; refer to \cite[Corollary~3.15(7)]{DMM1} and its proof.
Thus, $c_{\mu} \cup c_1^m$ is isotopic in $T_{-3, 2}(m)$
to the union of the core of the filled solid torus
and a regular fiber in $S^3 -\mathrm{int}N( T_{-3, 2} )$.
Thus, the lens space $T_{-3, 2}(m)$ has
a non-degenerate Seifert fibration
with $c_{\mu}$ and $c_1^m$ regular fibers. 
As shown in Figure~\ref{T-32m-5-7}, 
$c_{\mu} \cup c_1^{-5}$ is the $(2, -4)$ torus link and 
cobounds an annulus.
Note that such an annulus
cannot be disjoint from $T_{-3, 2}$ 
because $|\mathrm{lk}(c_{\mu}, T_{-3, 2})| = 1 \ne |m+ 1| 
= |\mathrm{lk}(c_1^m, T_{-3, 2})|$ 
for $m \in \{-5, -7\}$ by Proposition~\ref{irrelevant annular}(1).  
Thus $\{c_{\mu}, c_1^{-5}\}$ is an annular pair. 
On the other hand, 
since $c_{\mu} \cup c_1^{-7}$ is the Whitehead link, 
$\{c_{\mu}, c_1^{-7}\}$ is not an annular pair. \par

(2) The lens space $T_{-3, 2}(m)$ has a Seifert fibration
with $s_{-3}$ an exceptional fiber, 
so that 
$V' = T_{-3, 2}(m) - \mathrm{int}N(s_{-3})$ is a solid torus
and obtained from the solid torus
$V = S^3 - \mathrm{int}N(s_{-3})$ by $m$--surgery on $T_{-3, 2}$.
Since $T_{-3, 2}$ is the $(-3, 2)$ cable of $V$,
a meridian of $V'$ represents $[\mu] \pm 2(2[\lambda] -3[\mu])
= \pm ( 4[\lambda] + m[\mu] ) \in H_1( \partial V )$,
where $(\mu, \lambda)$ is a preferred meridian--longitude pair
of $V (\subset S^3)$.
Note that $s_{-3} \cup c_2^m$ 
in Figure~\ref{T-32m-5-7_2}  is isotopic in $T_{-3,2}(m)$ to
the union of $s_{-3}$ and the $(p, 1)$ cable of $s_{-3}$ in $S^3$.
Since $(p, 1) \ne (4, m)$, 
$T_{-3, 2}(m)$ has a non-degenerate Seifert fibration in which
$s_{-3}$ and $c_2^m$ are fibers. 
As shown in Figure~\ref{T-32m-5-7_2}, 
$s_{-3} \cup c_2^m$ is the $2$--bridge link associated to 
$\frac{6p+4}{2p+1}$. 
It is a torus link exactly when $p = 0, -1$. 
Hence $s_{-3}$ and $c_2^m$ cobound
an annulus exactly when $p = 0, -1$. 
Since $|\mathrm{lk}(s_{-3}, T_{-3, 2})| = 2 \ne |m+2| 
= |\mathrm{lk}(c_2^m, T_{-3, 2})|$ for 
$m \in \{-5, -7\}$, 
the annulus must intersect $T_{-3, 2}$
by Proposition~\ref{irrelevant annular}. 
Hence, 
$\{s_{-3}, c_2^m\}$ is an annular pair if and only if $p = 0, -1$. 

(3) As in (2), 
$V' = T_{-3, 2}(m) - \mathrm{int}N(s_2)$ is a solid torus and
obtained from the solid torus $S^3 - \mathrm{int}N(s_2)$
by $m$--surgery on $T_{-3, 2}$.
Since $T_{-3, 2}$ is the $(-2, 3)$ cable of $V$,
a meridian of $V'$ is a $(9, m)$ cable of $s_2$. 
Since $(p, 1)\ne (9, m)$,
$T_{-3, 2}(m)$ has a non-degenerate Seifert fibration 
with $s_2$ and $c_3^m$ fibers. 
Figure~\ref{T-32m-5-7_3} shows that
$s_2 \cup c_3^m$ is the $2$--bridge link associated to 
$\frac{6p+10}{2p+3}$. 
It is a torus link exactly when $p = -1, -2$. 
Since $|\mathrm{lk}(s_2, T_{-3, 2})| = 3 \ne |m+3|
= |\mathrm{lk}(c_3^m, T_{-3, 2})|$ for 
$m \in \{-5, -7\}$, 
$\{s_2, c_3^m\}$ is an annular pair if and only if $p = -1, -2$. 
\QED{Lemma~\ref{seifert pair example}}

Let us determine which annular pairs given in
Lemmas~\ref{annularT-32} and \ref{seifert pair example}
are basic annular pairs. 

\begin{proposition}
\label{annular pairs T-32}
Let $\{ \alpha_m, \beta_m \}$ be an annular pair of seiferters
in Figures~\ref{T-32annularpair1}--\ref{T-32annularpair4}
and \ref{T-32m-5-7}--\ref{T-32m-5-7_3}, 
and assume that
$\{ \alpha_m, \beta_m \}$ is a basic annular pair of seiferters for 
$T_{-3, 2}$.
Then $m= -3, -4, -5, -6$, and one of the following holds. 

\begin{enumerate}
\item
$\{ \alpha_m, \beta_m \} = \{ c_1^{-3}, c_2^{-3} \}$
in Figure~\ref{T-32annularpair4} if $m = -3$. 
\item
$\{ \alpha_m, \beta_m \} = \{ c_1^{-4}, c_2^{-4}\}$,
$\{  c_1^{-4}, c_3^{-4} \}$ or $\{ c_2^{-4}, c_3^{-4} \}$ 
in Figure~\ref{T-32annularpair4} if $m = -4$. 
\item
$\{ \alpha_m, \beta_m \} = \{ c_2^{-5}, c_3^{-5} \}$
in Figure~\ref{T-32annularpair4},
$\{ s_{-3}, c_2^{-5} \}$ with $p = -1$ in Figure~\ref{T-32m-5-7_2}, or
$\{ s_2, c_3^{-5} \}$ with $p=-2$ in Figure~\ref{T-32m-5-7_3} if $m = -5$. 
\item
$\{ \alpha_m, \beta_m \} = \{ s_{-3}, c_3^{-6} \}$
in Figure~\ref{T-32annularpair2} if $m = -6$. 
\end{enumerate}

Conversely, the pairs in $(1)$--$(4)$ are
basic annular pairs of seiferters for $T_{-3, 2}$. 
\end{proposition}

\noindent
\textit{Proof of Proposition \ref{annular pairs T-32}}.

\textit{Case 1. $\{ \alpha_m, \beta_m \}$ is
an annular pair of seiferters in
Figures~\ref{T-32annularpair1}--\ref{T-32annularpair4}}.

The linking numbers between any two of
$T_{-3, 2}, \alpha_m$, and $\beta_m$
are given in Table~\ref{table of linking numbers} below. 

\begin{table}[htbp]
\renewcommand{\arraystretch}{1.2}
\caption{Linking numbers of seiferters in
Figures~\ref{T-32annularpair1}--\ref{T-32annularpair4}}
\label{table of linking numbers}

\begin{tabular}{|c|c|r@{}r@{}r@{}r@{}r@{}}
\hline
lk & $T_{-3, 2}$\ &\ $c_{\mu}$\ \vline&\ $s_{-3}$\ \vline& \ $s_2$ \vline& $c_1^m$\ \ \ \vline&$c_2^m$\ \ \  \vline \\ 
\hline
$c_3^m$ & $m+3$ & $1$ \vline& $1$ \vline& $*$ \vline&\ $m+4$ \vline&\ $m+4$ \vline\\
\hline
$c_2^m$ & $m+2$ & $1$ \vline& $*$ \vline& $2$ \vline&\ $m+3$ \vline&  \\
\cline{1-6}
$c_1^m$\ & $m+1$ & $*$ \vline& $2$ \vline& $3$ \vline& &\\
\cline{1-5}
$s_2$   & $ 3$   & &   & & & \\
\cline{1-2}
$s_{-3}$ & $2$   & &   & & & \\
\cline{1-2}
$c_{\mu}$\ & $1$ & &    & & & \\
\cline{1-2}
\end{tabular}
\end{table}

\begin{figure}[h]
\begin{center}
\includegraphics[width=0.9\linewidth]{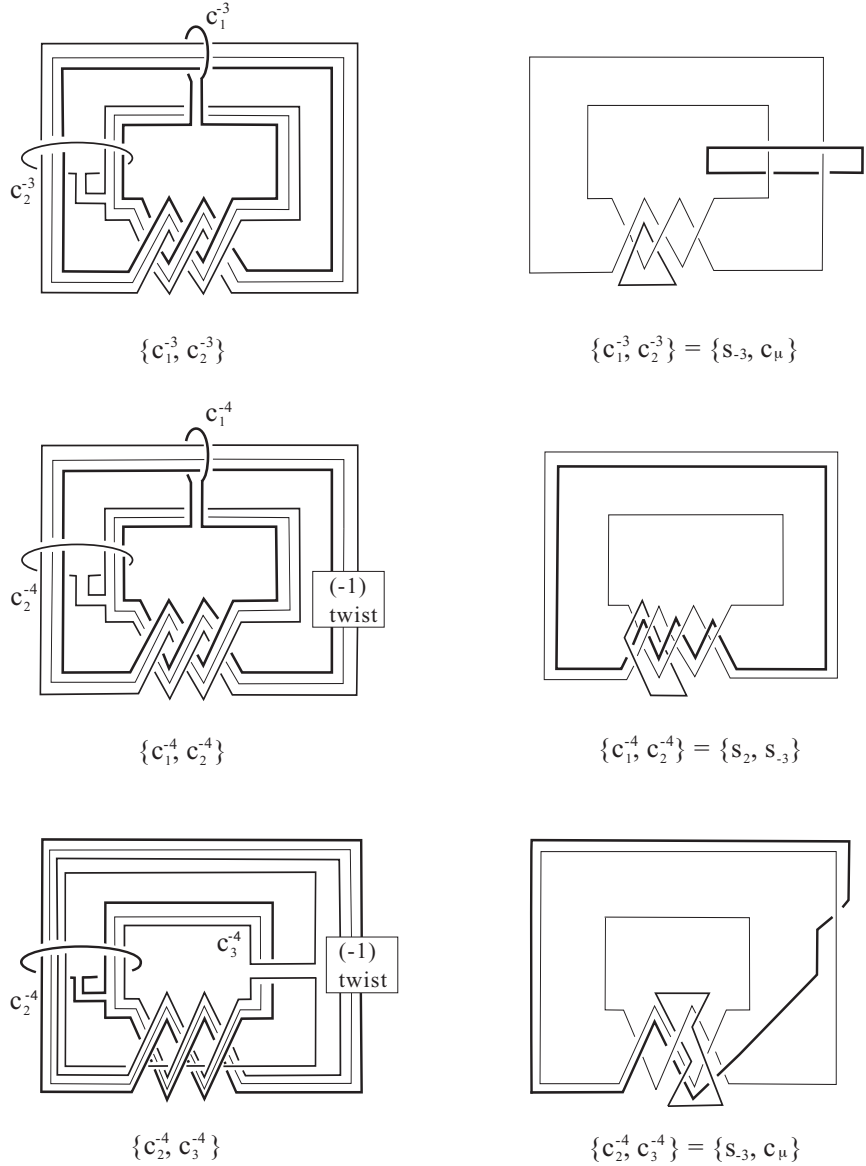}
\caption{$\{c_1^{-3}, c_2^{-3}\}$, $\{c_1^{-4}, c_2^{-4}\}$
and $\{c_2^{-4}, c_3^{-4}\}$ 
are basic annular pairs for $T_{-3,2}$.}
\label{T-32basicannular1}
\end{center}
\end{figure}

\begin{figure}[h]
\begin{center}
\includegraphics[width=0.9\linewidth]{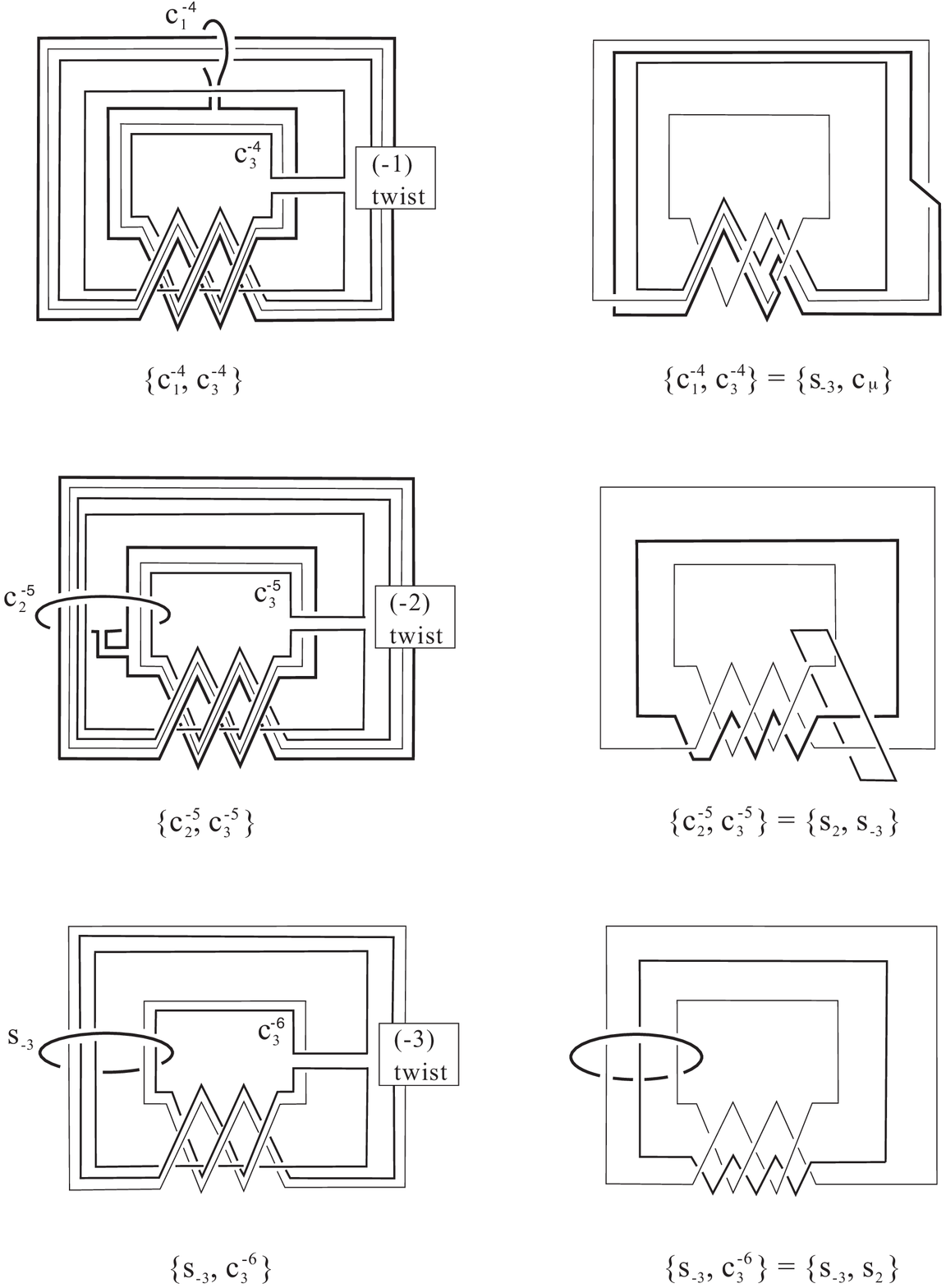}
\caption{$\{ c_1^{-4}, c_3^{-4} \}$, $\{c_2^{-5}, c_3^{-5} \}$,
and $\{ s_{-3}, c_3^{-6} \}$ 
are basic annular pairs for $T_{-3, 2}$.}
\label{T-32basicannular2}
\end{center}
\end{figure}

Since $\{ \alpha_m, \beta_m \}$
is a basic annular pair of seiferters for $T_{-3, 2}$, 
the triple 
$( |\mathrm{lk}(T_{-3, 2}, \alpha_m)|$,
$|\mathrm{lk}(T_{-3, 2}, \beta_m)|$, 
$|\mathrm{lk}(\alpha_m, \beta_m)| )$
is equal to  
$(1, 2, 0), (1, 3, 0)$ or $(2, 3, 1)$ according as
$(\alpha_m, \beta_m) =  (c_{\mu}, s_{-3})$, 
$(c_{\mu}, s_2)$
or 
$(s_{-3}, s_2)$. 
Checking the triples of linking numbers 
for the nine possible pairs,
we obtain all cases listed in
Proposition~\ref{annular pairs T-32}(1)--(4)
and one unexpected case $\{ s_{-3}, c_3^0 \}$.
However, the latter is not a basic annular pair
because $c_3^0( =c^2)$ is a hyperbolic seiferter for $(T_{-3, 2}, 0)$
by Proposition~\ref{cm is nonbasic}(1).
Conversely,
the pairs in Proposition~\ref{annular pairs T-32}(1)--(4)
are basic annular pairs of seiferters,
as shown in
Figures~\ref{T-32basicannular1} and \ref{T-32basicannular2}. 

\textit{Case 2. $\{ \alpha_m, \beta_m \}$ is
an annular pair of seiferters in
Figures~\ref{T-32m-5-7}--\ref{T-32m-5-7_3}}.

Checking whether $c_1^m, c_2^m, c_3^m$ $(m = -5, -7)$ are basic seiferters
by Table~\ref{table of c^m}, 
we see that $m = -5$, and
$\{ c_{\mu}, c_1^m  \}$ in Figure~\ref{T-32m-5-7} is not
a basic annular pair of seiferters.
We also see that $\{ \alpha_m, \beta_m \}$ is either
$\{ s_{-3}, c_2^m  \}$ in Figure~\ref{T-32m-5-7_2} or
$\{ s_2, c_3^m \}$ in Figure~\ref{T-32m-5-7_3}, and in either case
it is the basic annular pair $\{s_{-3}, s_2\}$.
If $\{ s_{-3}, c_2^{-5} \}$ in Figure~\ref{T-32m-5-7_2} 
is $\{s_{-3}, s_2 \}$,
then $p\in \{0, -1\}$ by Lemma~\ref{seifert pair example}
and the 2--bridge link $s_{-3} \cup c_2^{-5}$
in the right-most figure of Figure~\ref{T-32m-5-7} is a Hopf link.
It follows $p = -1$.
Similarly, 
if $\{ s_2,  c_3^{-5} \}$ in Figure~\ref{T-32m-5-7_3} 
is $\{s_{-3}, s_2\}$,
then we see $p = -2$
from the right-most figure of Figure~\ref{T-32m-5-7}.
Proposition~\ref{annular pairs T-32} follows from 
Lemma~\ref{lem:Hopf links} below.
\QED{Proposition~\ref{annular pairs T-32}}

\begin{lemma}
\label{lem:Hopf links}
\begin{enumerate}
\item $\{ s_{-3}, c_2^{-5} \}$ with $p =-1$ in Figure~\ref{T-32m-5-7}
is the basic annular pair $\{ s_{-3}, s_2 \}$.
\item $\{ s_2, c_3^{-5} \}$ with $p =-2$ in Figure~\ref{T-32m-5-7_3}
is the basic annular pair $\{ s_{-3}, s_2 \}$.
\end{enumerate}
\end{lemma}

\noindent
\textit{Proof of Lemma~\ref{lem:Hopf links}}.  
Although we can depict the isotopies showing the lemma, 
we give a proof using a more general argument.

(1) We show that the exterior
$X = S^3 -\mathrm{int}N( s_{-3} \cup c_2^{-5} \cup T_{-3, 2})$ is homeomorphic
to $(\mbox{the twice punctured disk})\times S^1$.
Then, by \cite{BuMu} $s_{-3} \cup c_2^{-5} \cup T_{-3, 2}$ is a union of
fibers of some Seifert fibration of $S^3$.
This implies the desired result.

Recall that
$c_2^m$ is isotopic in $S^3 -\mathrm{int}N(s_{-3} \cup T_{-3, 2})$
to a band sum of a simple closed curve
in $\partial N(s_{-3})$ with slope $p$
and one in $\partial N(T_{-3, 2})$ with slope $m= -5$;
let $b$ be the band used in this band sum,
where $b \subset S^3 -\mathrm{int}N(s_{-3} \cup T_{-3,2})$.
Then $s_{-3}, c_2^m$ and $T_{-3, 2}$
cobound an obvious planar surface.
By restricting this surface in
$X = S^3 -\mathrm{int}N(s_{-3} \cup c_2^m \cup T_{-3, 2})$,
we obtain a twice punctured disk $S$ properly embedded in $X$.
The boundary slopes of $S$ are $p = -1$ in $\partial N(s_{-3})$,
$m = -5$ in $\partial N(T_{-3, 2})$,
and $m + p + 2\mathrm{lk}(s_{-3}, T_{-3, 2}) = -2$ in $\partial N(c_2^{-5})$
by Proposition~\ref{move}.
Take a collar neighborhood $S\times I$ of $S$ in $X$
such that $S\times\{ 0 \} =S$ and 
$\partial S \times I = \partial X \cap S \times I$.
Let $Y$ be the closure of $X - S\times I$.
Then $Y$ is homeomorphic to the closure of
$S^3 - N(s_{-3} \cup T_{-3,2}) - N(a)$,
where $a$ is the core of the band $b$.
The arc $a$ is a non-separating arc in an essential annulus in
the cable space $S^3 -\mathrm{int}N(s_{-3} \cup T_{-3, 2})$
which splits the cable space into a solid torus.
This implies that $Y$ is a handlebody of genus 2.
We aim to prove $(Y, S) \cong ( S \times I, S \times \{ 0 \} )$.
This implies that $X$ is a fiber bundle with $S$ a fiber.
Since any self-homeomorphism of the twice punctured disk
with its boundary setwise invariant is isotopic to the identity,
we see $X \cong S \times S^1$.

For simplicity, set $K_1 = s_{-3}$, $K_2 = c_2^{-5}$, 
$K_3 = T_{-3, 2}$, and $\gamma_i = S \cap \partial N(K_i)$,
where $i=1,2,3$.
We also denote by $[\gamma_i]$ the slope of $\gamma_i$ in
$\partial N(K_i)$.
For a set $\mathcal{C}$ of disjoint simple closed curves on
a 3-manifold $M$, $\tau(M; \mathcal{C})$ denote
$M$ with 2--handles added along the loops in $\mathcal{C}$.

\begin{claim}
\label{cl:handle addition}
For any proper subset $\mathcal{C}'$ of 
$\{ \gamma_1, \gamma_2, \gamma_3 \}$,
$\tau(Y; \mathcal{C}')$ is a handlebody.
\end{claim}

By \cite[Theorem~2]{Go87} Claim~\ref{cl:handle addition} implies
$(Y, S) \cong ( S \times I, S \times \{ 0 \} )$ as desired.
We prove Claim~\ref{cl:handle addition} by relating
$\tau(Y; \{ \gamma_i \})$ with the Dehn surgery $K_i([ \gamma_i ])$.

\begin{claim}
\label{cl:handle addition2}
$\tau(Y, \{ \gamma_i \}) \cong K_i( [\gamma_i] ) -\mathrm{int}N( K_j )$,
where $i\ne j$.
\end{claim}

\noindent
\textit{Proof of Claim~\ref{cl:handle addition2}.}
Let $V$ be the filled solid torus in $K_i( [\gamma_i] )$.
Cut $V$ by disks bounded by two meridians
$S \times \{0, 1\} \cap N(K_i)$ into two 3--balls $W_1, W_2$.
We may assume $W_1 \cap Y = \partial N(K_i) \cap Y$,
so that $W_1$ (resp.\ $W_2$) is attached to $Y$
(resp.\ $S\times I$) as a 2--handle.
Then, for $\{i, \alpha, \beta\} =\{1, 2, 3\}$,
$K_i( [\gamma_i] )
= X \cup V \cup N(K_{\alpha}) \cup N(K_{\beta})
= ( Y \cup W_1) \cup ( S\times I \cup W_2)
\cup N(K_{\alpha}) \cup N(K_{\beta})$.
Note $Y\cup W_1 \cong \tau( Y; \{ \gamma_i \} )$. 
In $K_i( [\gamma_i] )$, $S$ with $\gamma_i$ capped off
by a meridian disk of $V$ is an annulus connecting
$\gamma_{\alpha} (\subset N(K_{\alpha}) )$ and
$\gamma_{\beta}( \subset N(K_{\beta}) )$.
Hence $( S\times I \cup W_2) \cup N(K_{\alpha}) \cup N(K_{\beta})$
is a solid torus isotopic in $K_i( [\gamma_i] )$ to
both $N(K_{\alpha})$ and $N(K_{\beta})$.
We then see that
$\tau( Y; \{ \gamma_i \} ) \cong K_i( [\gamma_i] ) -\mathrm{int}N(K_j)$,
where $j \in \{ \alpha, \beta \}$.
\QED{Claim~\ref{cl:handle addition2}}

Since $s_{-3} \cup c_2^{-5}$ is a Hopf link,
$s_{-3}( [\gamma_1] ) -\mathrm{int}N(c_2^{-5})$ and
$c_2^{-5}( [\gamma_2] ) -\mathrm{int}N(s_{-3})$ are solid tori.
The manifold $T_{-3, 2}( [\gamma_3] ) -\mathrm{int}N(s_{-3})$ is
obtained from the cable space $S^3 -\mathrm{int}N( s_{-3} \cup T_{-3,2})$
 by Dehn-filling $\partial N(T_{-3, 2})$ along $\gamma_3$.
Since $[\gamma_3] = -5$, $\gamma_3$ in $\partial N(T_{-3, 2})$
meets a fiber of a Seifert fibration of the cable space exactly once. 
Hence,
$T_{-3, 2}( [\gamma_3] ) -\mathrm{int}N(s_{-3})$ is a solid torus,
so that $\tau( Y, \{ \gamma_i \})$ is a solid torus for any $i$.

\begin{claim}
\label{cl:handle addition3}
For a pair $\{\alpha, \beta \} \subset \{1, 2, 3\}$,
let $M_{\alpha \beta}$ be the manifold obtained from $S^3$
by applying $[\gamma_{\alpha}]$--surgery on $K_{\alpha}$ and 
$[\gamma_{\beta}]$--surgery on $K_{\beta}$.
Then, $\tau(Y; \{ \gamma_{\alpha}, \gamma_{\beta} \})$ is homeomorphic to
a punctured $M_{\alpha \beta}$.
\end{claim}

\noindent
\textit{Proof of Claim \ref{cl:handle addition3}.}
As in the proof of Claim~\ref{cl:handle addition2},
the filled solid tori in $M_{\alpha \beta}$ are decomposed into 
$3$--balls which are 2--handles attached to
$Y$ and $S\times I$.
Then, we see that
$M_{\alpha \beta} = M_1 \cup M_2 \cup N(K_i)$,
where $i \not \in \{\alpha, \beta\}$,
$M_1$ (resp.\ $M_2$) is
$Y$ (resp.\ $S\times I$) with two 2--handles added along
the two annuli $Y \cap \partial N(K_{\alpha} \cup K_{\beta})$
(resp.\ $S\times I \cap \partial N(K_{\alpha} \cup K_{\beta})$).
Note $M_1 \cong \tau( Y; \{ \gamma_{\alpha}, \gamma_{\beta} \})$ and
$M_2 \cong \tau( S \times I; \{ \gamma_{\alpha}, \gamma_{\beta} \})$.
Since $M_2$ is a 3--ball and attached to $N(K_i)$ along $\gamma_i$
as a 2--handle,
$M_2 \cup N(K_i)$ is a 3--ball.
Hence Claim~\ref{cl:handle addition3} follows.
\QED{Claim~\ref{cl:handle addition3}}

Since $N( K_{\alpha} )$ is isotopic in $K_{i}( [\gamma_i] )$
to $N( K_{\beta} )$
with $\gamma_{\alpha}$ sent to $\gamma_{\beta}$,
we see $M_{\alpha \beta} \cong M_{12}$.
Since $s_{-3} \cup c_2^{-5}$ is a Hopf link,
$M_{12}$, which is the result of
$(-1)$--surgery on $s_{-3}$ and $(-2)$--surgery on $c_2^{-5}$,
is the 3--sphere.
It then follows from Claim~\ref{cl:handle addition3} that
$\tau(Y; \{ \gamma_{\alpha}, \gamma_{\beta} \})$ is a 3--ball.
This completes the proof of Claim~\ref{cl:handle addition}. 
\QED{Claim~\ref{cl:handle addition}}

(2) The arguments in (1) apply after some replacement.
There is a twice punctured disk $S$ properly embedded in
$X = S^3 -\mathrm{int}N(s_2 \cup c_3^{-5} \cup T_{-3, 2})$
such that the boundary slopes of $S$ are
$p= -2$ in $\partial N(s_2)$, $m =-5$ in $\partial N(T_{-3,2})$,
and $m + p + 2\mathrm{lk}(s_2, T_{-3, 2}) = -1$ in $N(c_3^{-5})$.
Note also that $s_2 \cup c_3^{-5}$ is a Hopf link.
The arguments after Claim~\ref{cl:handle addition}
hold with $s_{-3}$ and $c_2^{-5}$ replaced
by $c_3^{-5}$ and $s_2$, respectively.
\QED{Lemma~\ref{lem:Hopf links}}

Theorem~3.24 in \cite{DMM1} shows that
under some conditions an annular pair of seiferters 
is either hyperbolic or basic.
Using this theorem and Proposition~\ref{annular pairs T-32}, 
we give a sufficient condition for annular pairs in
Figures~\ref{T-32annularpair1}--\ref{T-32annularpair4}
to be hyperbolic.

\begin{corollary}
\label{hyperbolic annular pairs T-32}
\begin{enumerate}
\item
Let $m$ be an integer other than $-5, -6, -7$, 
Then, an annular pair of seiferters 
in Figures~\ref{T-32annularpair1}--\ref{T-32annularpair4}
for $(T_{-3, 2}, m)$ is hyperbolic 
if and only if it is not 
$\{ c_1^{-3}, c_2^{-3} \}$ in Figure~\ref{T-32annularpair4} 
with $m= -3$ or
$\{ c_1^{-4}, c_2^{-4} \}$, $\{ c_1^{-4}, c_3^{-4} \}$, 
$\{ c_2^{-4}, c_3^{-4} \}$ 
in Figure~\ref{T-32annularpair4} with $m = -4$. 
\item
If $m \le -8$ or $-1 \le m$, then
nine annular pairs of seiferters for $(T_{-3, 2}, m)$
in Figures~\ref{T-32annularpair1}--\ref{T-32annularpair4} are
all hyperbolic.
\end{enumerate}
\end{corollary}

\noindent
\textit{Proof of Corollary~\ref{hyperbolic annular pairs T-32}}.
(1) Let $\{ \alpha_m, \beta_m \}$ be 
an annular pair of seiferters for $(T_{-3, 2}, m)$ in
Figures~\ref{T-32annularpair1}--\ref{T-32annularpair4}.
Since $m \ne -6$, $T_{-3, 2}(m)$ is not
a connected sum of lens spaces,
so that $(T_{-3, 2}, m)$ is a Seifert fibered surgery.
The assumption $m \not\in \{-5, -7\}$ implies that
$T_{-3, 2}(m)$ is not a lens space.
Since $\alpha_m\cup \beta_m$ is isotopic in $T_{-3, 2}(m)$ to
a basic annular pair for $T_{-3,2}$,
$\alpha_m$ and $\beta_m$ are exceptional fibers 
in a Seiferter fibration of $T_{-3, 2}(m)$.
Then, Theorem~3.24 in \cite{DMM1} shows that
$\{ \alpha_m, \beta_m \}$ is either hyperbolic or basic.
Now assertion~(1) follows from Proposition~\ref{annular pairs T-32}.

(2) Since $m \le -8$ or $-1 \le m$,
we can apply assertion~(1).
Then, since $m \not\in \{-3, -4\}$, the desired result follows.
\QED{Corollary~\ref{hyperbolic annular pairs T-32}}

\begin{remark}
\label{annular(T-32, -4)}
Corollary~\ref{hyperbolic annular pairs T-32}
shows that $(T_{-3, 2}, -4)$ has six hyperbolic,
annular pairs of seiferters
$\{ c_{\mu}, c_2^{-4} \}$, $\{ c_{\mu}, c_3^{-4} \}$, 
$\{ s_{-3}, c_1^{-4} \}$, $\{ s_{-3}, c_3^{-4} \}$, 
$\{ s_2, c_1^{-4} \}$, $\{ s_2, c_2^{-4} \}$.
Note that each of these consists of
basic seiferters for $T_{-3, 2}$
by Proposition~\ref{cm is nonbasic}. 
\end{remark}

\section{Strongly invertible knots that do not arise from
the primitive/Seifert--fibered construction}
\label{section:Song}

We first review the definition of primitive/Seifert-fibered construction introduced by Dean \cite{D2}.
Let $K$ be a knot in a genus 2 Heegaard surface $F$ of $S^3 =H \cup_F H'$,
where $H, H'$ are genus 2 handlebodies.
The surface slope $\gamma_{K, F}$ of $K (\subset F)$
is the isotopy class in
$\partial N(K)$ represented by a component of $\partial N(K) \cap F$
which is parallel to $K$.
We denote by $H[ K ]$ (resp.\ $H'[ K ]$)
the 3--manifold $H$ (resp.\ $H'$) with a 2--handle added along $K$.
Note that the surgered manifold $K(\gamma_{K, F})$
is the union of $H[ K ]$ and $H'[ K ]$.
The knot $K$ is said to be \textit{primitive/Seifert-fibered} with respect
to $F$ if $H[ K ]$ is a solid torus
and $H'[ K ]$ is a Seifert fiber space over $D^2(p, q)$, where
$p, q \ge 2$.
If $K$ is primitive/Seifert-fibered, then
$K(\gamma_{K, F})$ is a lens space, a small Seifert fiber space,
or a connected sum of two lens spaces;
$(K, \gamma_{K, F})$ is a Seifert surgery.
We say that a Seifert surgery $(K, m)$ arises from
the \textit{primitive/Seifert-fibered construction} 
if $K$ is isotopic to a knot $L$ in a genus 2 Heegaard surface $F$
of $S^3$ in such a way that $L$ is primitive/Seifert-fibered
with respect to $F$ and the surface slope $\gamma_{L, F}$
coincides with $m$.
This construction of Seifert surgeries is a modification of
Berge's primitive/primitive construction \cite{Berge2}
of lens space surgeries.
Primitive/primitive knots and
primitive/Seifert-fibered knots have tunnel number one,
and thus are strongly invertible \cite{D2}.

Although all known lens space surgeries arise from
primitive/primitive constructions,
there are infinite families of small Seifert fibered surgeries
none of which arises from
the primitive/Seifert-fibered construction \cite{MMM, DMM1, Tera}.
The simplest example is the $1$--surgery on 
the $(-3, 3, 5)$ pretzel knot \cite{MMM}. 
For any $(K, m)$ in the families found in \cite{MMM, DMM1, Tera}, 
$K$ is not strongly invertible. 
It is natural to raise the following question,
and Song \cite{Song} gives an example.

\begin{question}
\label{invertiblePS}
Does there exist a small Seifert fibered surgery
on a strongly invertible knot
which does not arise from the primitive/Seifert-fibered construction?
\end{question}

\begin{example}[{\cite{Song}}]
\label{Song's example}
$(-1)$--surgery on the strongly invertible pretzel knot $P(3, -3, -3)$ 
is a small Seifert fibered surgery which does not arise from 
the primitive/Seifert-fibered construction.
\end{example}

By twisting $(P(3,-3,-3), -1)$ along a seiferter,
we extend Song's example to a one-parameter family of Seifert surgeries
which give an affirmative answer to Question~\ref{invertiblePS}.
We first show that Song's example is obtained by
twisting a Seifert surgery on a trefoil knot.

\begin{figure}[htbp]
\begin{center}
\includegraphics[width=0.28\linewidth]{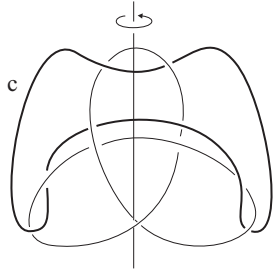}
\caption{Seiferter $c$ for $(T_{-3, 2}, -1)$.}
\label{nonpmseiferter}
\end{center}
\end{figure}

\begin{lemma}
\label{P(3,-3,-3)}
Let $c$ be the trivial knot given in Figure~\ref{nonpmseiferter}.
Then the following hold.
\begin{enumerate}

\item
The trivial knot $c$ is a hyperbolic seiferter for 
the Seifert surgery $(T_{-3, 2}, -1)$.

\item
$(-1)$--twist along $c$ converts $(T_{-3, 2}, -1)$ to 
the Seifert surgery $(P(3, -3, -3), -1)$. 

\item
The seiferter $c$ cannot be obtained from any basic seiferter for 
$T_{-3, 2}$ by a single $(-1)$--move. 

\end{enumerate}
\end{lemma}

\begin{figure}[htbp]
\begin{center}
\includegraphics[width=0.4\linewidth]{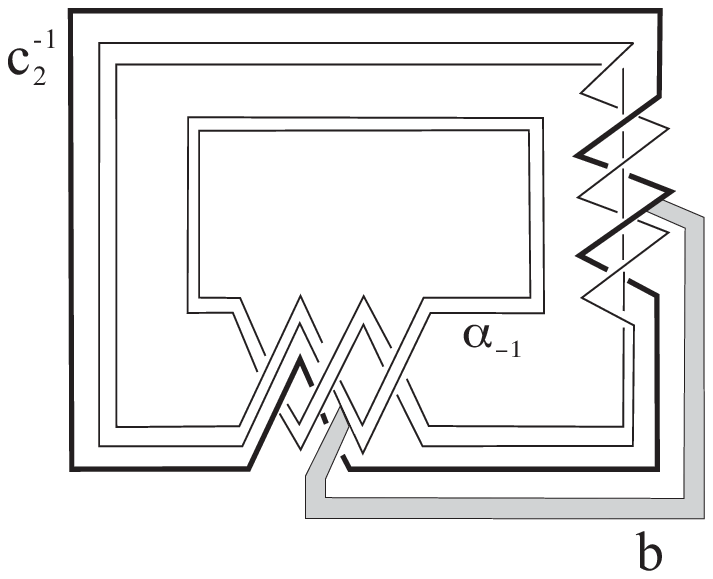}
\caption{}
\label{c0band}
\end{center}
\end{figure}

\begin{figure}[htbp]
\begin{center}
\includegraphics[width=0.9\linewidth]{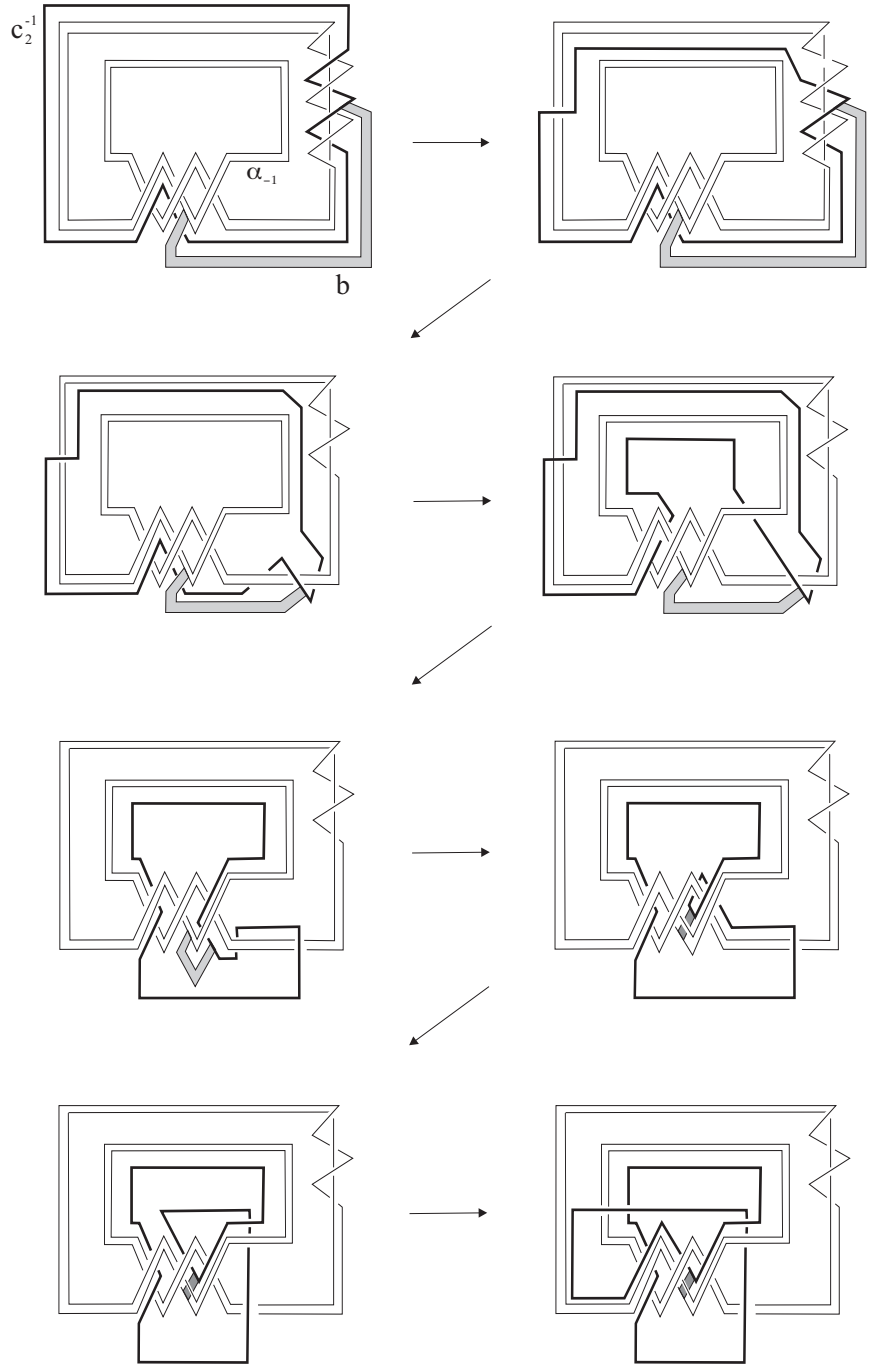}
\caption{Isotopy of $T_{-3, 2} \cup c_2^{-1}$}
\label{tp12A}
\end{center}
\end{figure}

\begin{figure}[htbp]
\begin{center}
\includegraphics[width=0.9\linewidth]{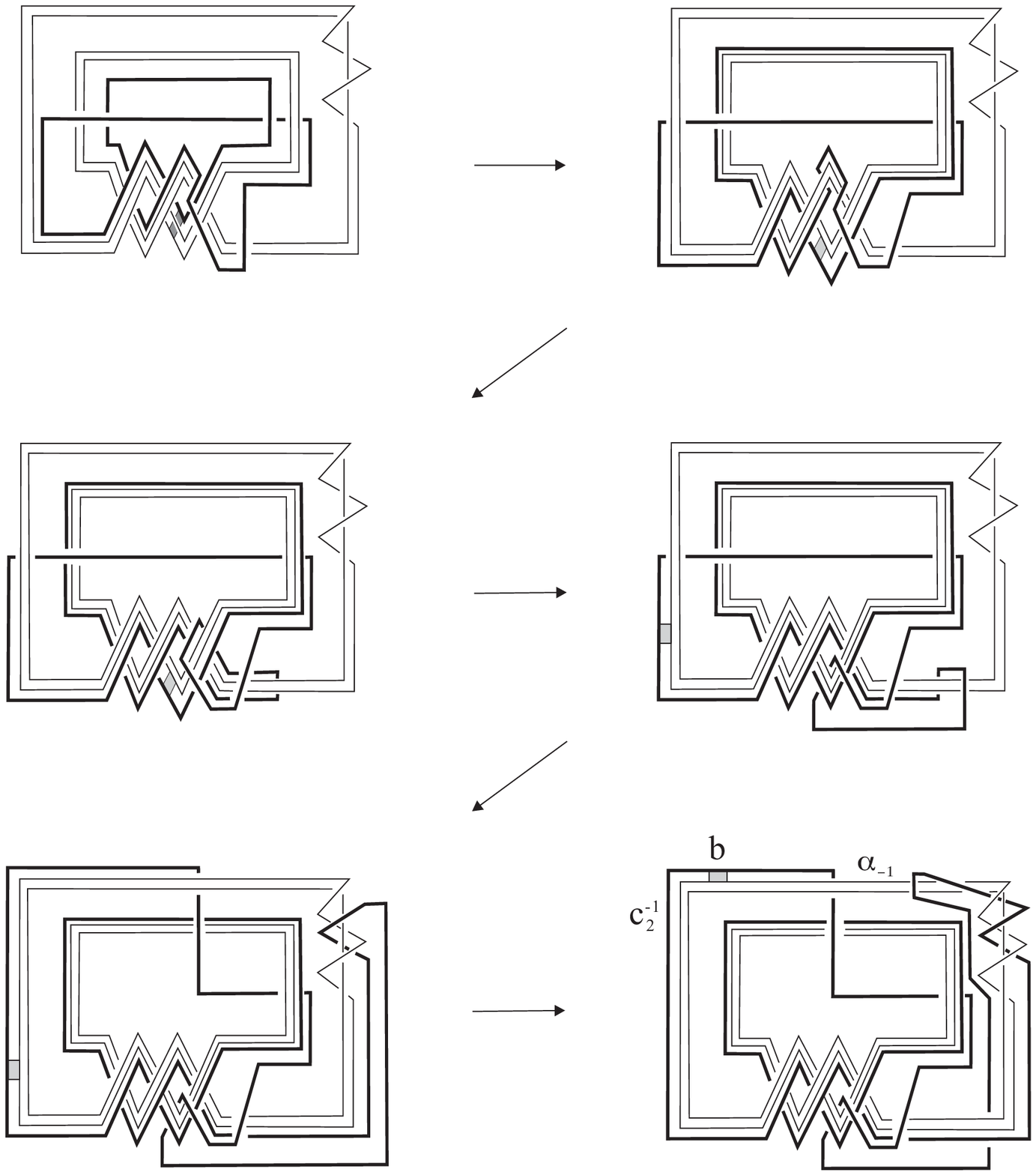}
\caption{Isotopy of $T_{-3, 2} \cup c_2^{-1}$; continued from Figure~\ref{tp12A}}
\label{tp12B}
\end{center}
\end{figure}

\begin{figure}[htbp]
\begin{center}
\includegraphics[width=0.9\linewidth]{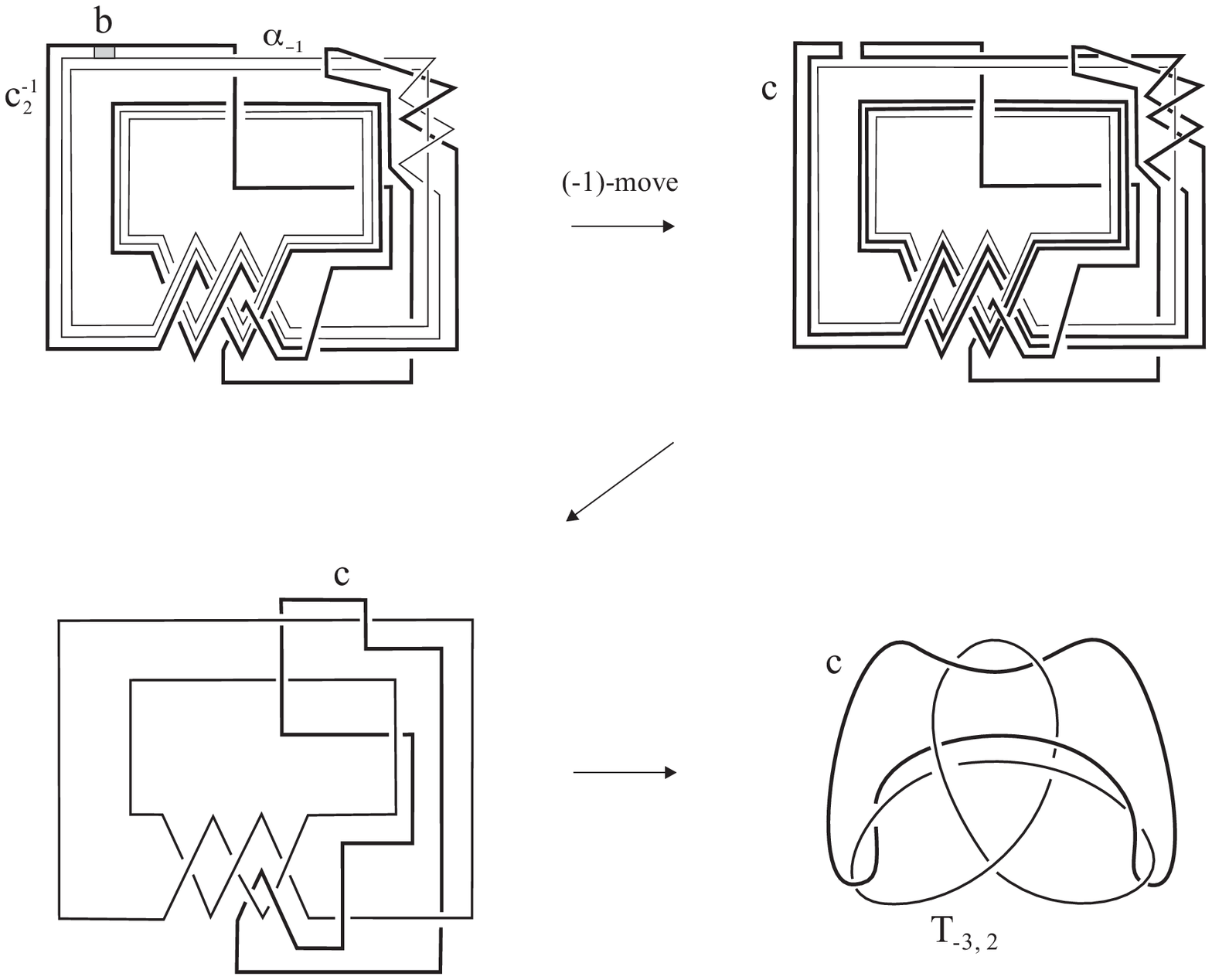}
\caption{}
\label{tpmove}
\end{center}
\end{figure}

\noindent
\textit{Proof of Lemma~\ref{P(3,-3,-3)}.}
(1) As shown in Figure~\ref{tc2},
the trivial knot $c_2^{-1}$ in Figure~\ref{c0band}
is a seiferter for $(T_{-3, 2}, -1)$, 
and an exceptional fiber of index $3$ in $T_{-3, 2}(-1)$. 
Let us take the band $b$ as in Figure~\ref{c0band},
which connects $c_2^{-1}$ and a simple closed curve $\alpha_{-1}$ 
on $\partial N(T_{-3, 2})$ with slope $-1$.
Isotope $c_2^{-1} \cup b$ as described in Figures~\ref{tp12A}, \ref{tp12B}.
Then Figure~\ref{tpmove} shows that the $(-1)$--move via the band $b$
converts $T_{-3, 2} \cup c_2^{-1}$ to $T_{-3, 2} \cup c$.
It follows that $c$ is a seiferter for $(T_{-3, 2}, -1)$
and an exceptional fiber of index $3$ in $T_{-3, 2}(-1)$.
Then, by \cite[Corollary~3.15]{DMM1} 
$c$ is a hyperbolic seiferter or a basic seiferter
for $(T_{-3, 2}, -1)$.
However, $c$ is not a basic seiferter for $T_{-3, 2}$
because $\mathrm{lk}(T_{-3, 2}, c) =0$.
This establishes assertion~(1). 
Note that the seiferter $c$ is obtained
from the basic seiferter $s_{-3}$ by applying $(-1)$--moves twice. 

\begin{figure}[htbp]
\begin{center}
\includegraphics[width=0.9\linewidth]{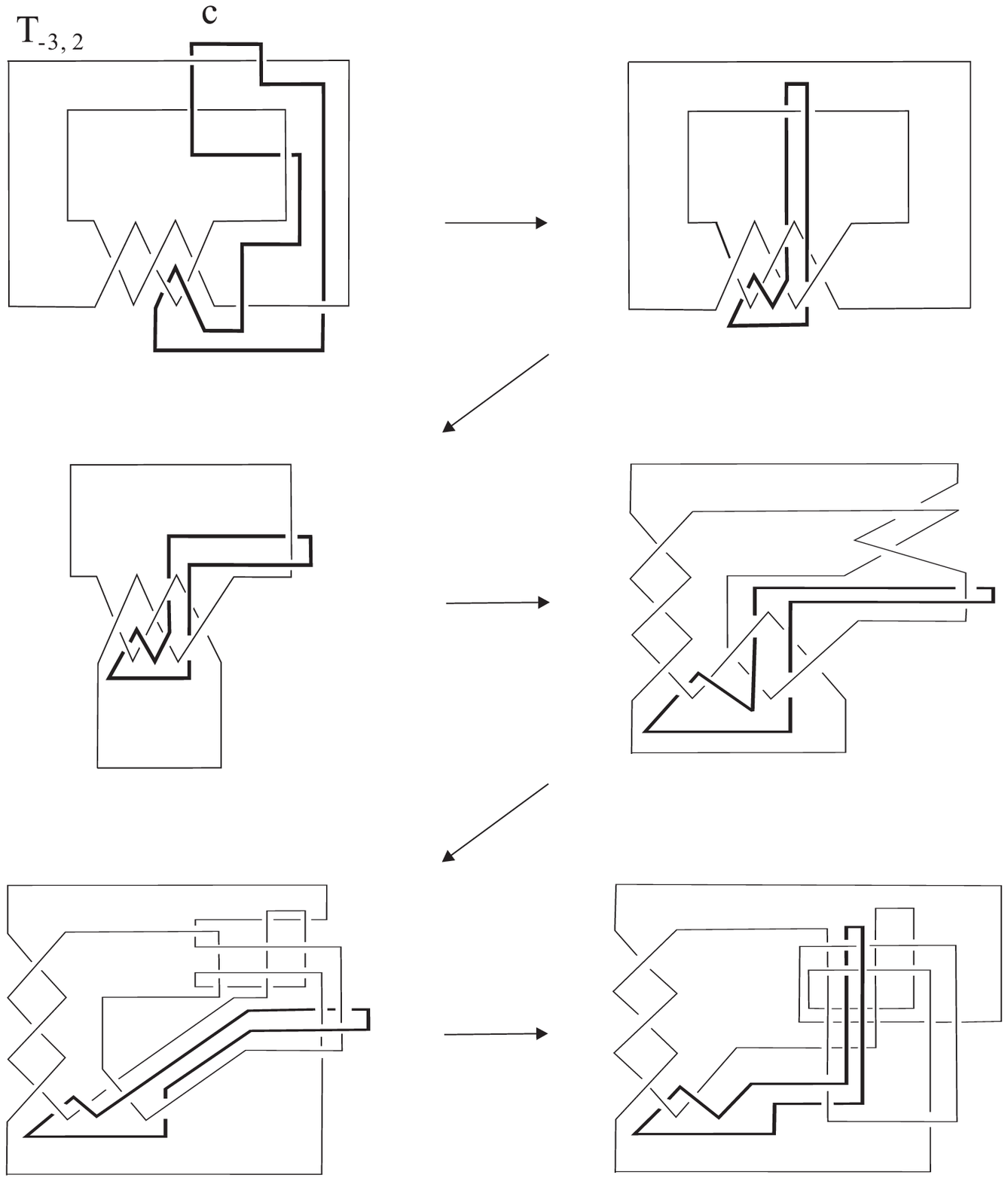}
\caption{Isotopy of $T_{-3, 2} \cup c$}
\label{tpisotopy}
\end{center}
\end{figure}

\begin{figure}[htbp]
\begin{center}
\includegraphics[width=0.9\linewidth]{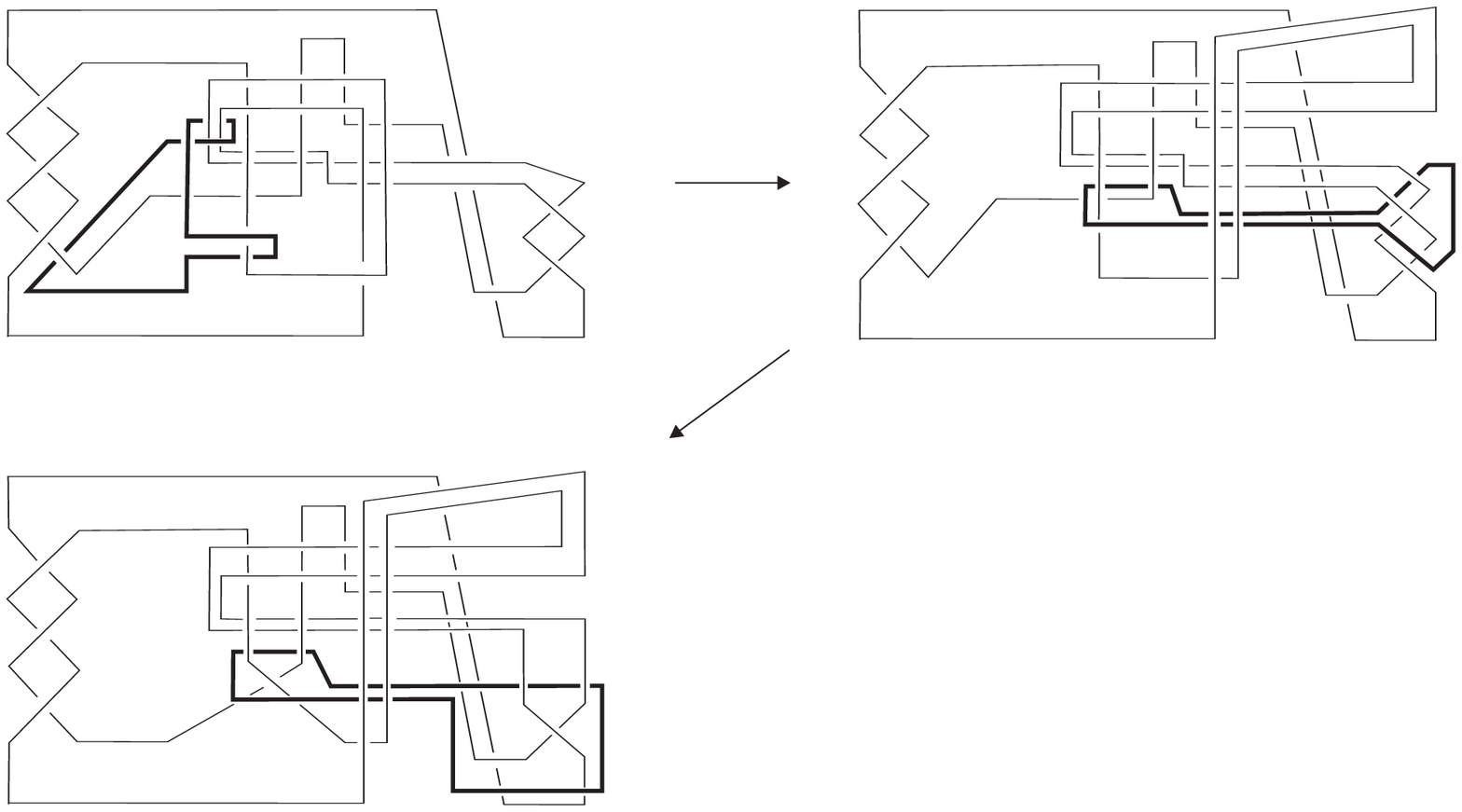}
\caption{Isotopy of $T_{-3, 2} \cup c$; 
continued from Figure~\ref{tpisotopy}}
\label{tpisotopy2}
\end{center}
\end{figure}

\begin{figure}[htbp]
\begin{center}
\includegraphics[width=0.9\linewidth]{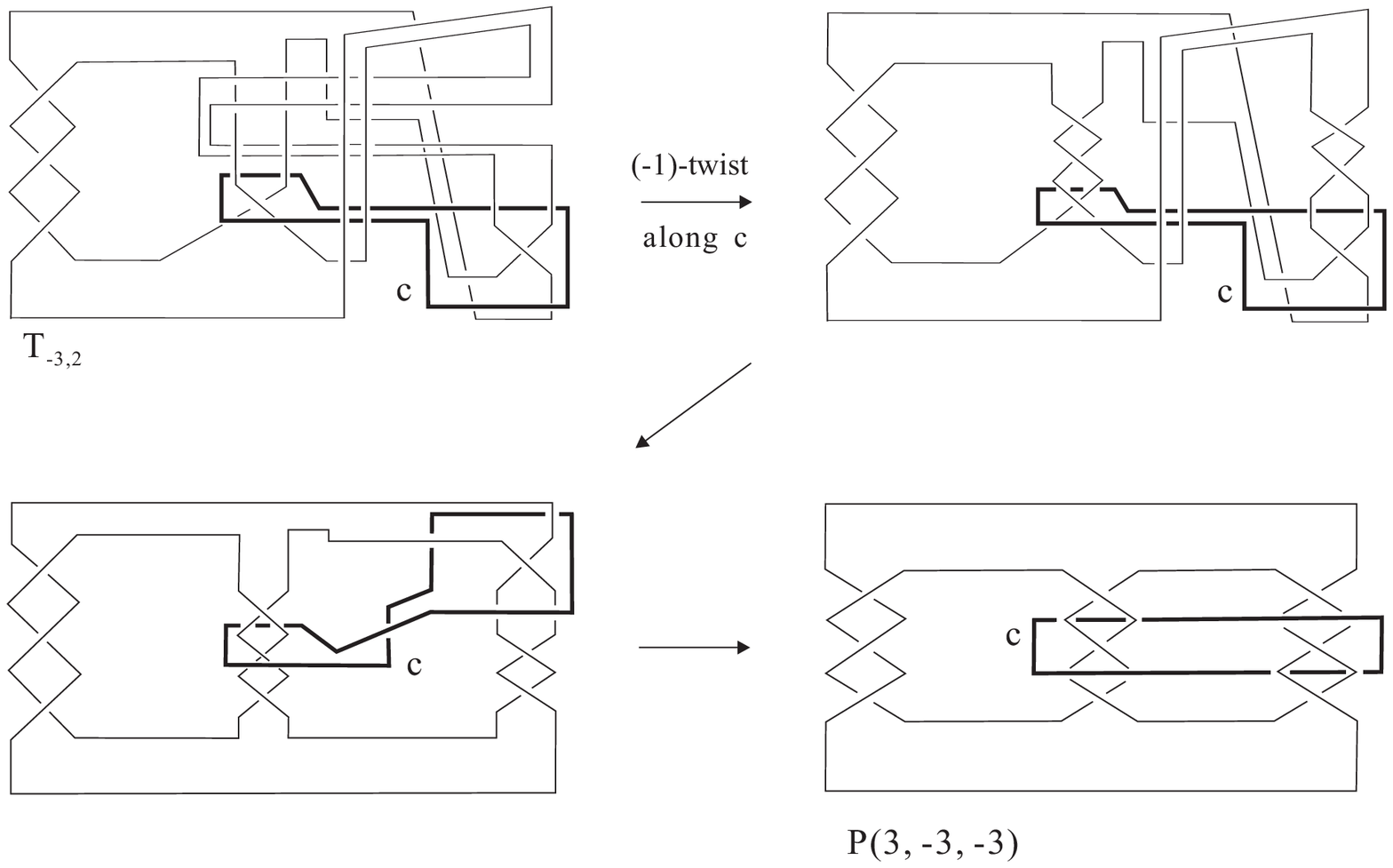}
\caption{$(-1)$--twist along $c$ converts $T_{-3, 2}$ to $P(3, -3, -3)$.}
\label{tptwist}
\end{center}
\end{figure}

(2) Isotope the link $T_{-3, 2} \cup c$ as in 
Figures~\ref{tpisotopy} and \ref{tpisotopy2},
and twist $T_{-3, 2}$ $(-1)$--time along the seiferter $c$.
We then obtain $P(3, -3, -3)$ as required; 
see Figure~\ref{tptwist}.
Since $\mathrm{lk}(T_{-3, 2}, c) =0$,
twists on $(T_{-3, 2}, -1)$ along $c$ do not change
the surgery coefficient. 
Hence, 
$(-1)$--twist along $c$ converts $(T_{-3, 2}, -1)$ to 
$(P(3, -3, -3), -1)$. \par

(3) Recall that $|\mathrm{lk}(T_{-3, 2}, c_{\mu})| = 1$, 
$|\mathrm{lk}(T_{-3, 2}, s_{-3})| = 2$, and 
$|\mathrm{lk}(T_{-3, 2}, s_{2})| = 3$.
If $c$ were obtained from $c_{\mu}$, $s_{-3}$, or $s_2$
by a single $(-1)$--move,  
then Proposition~\ref{move}(2) 
with $\varepsilon = \pm 1$ and $m= -1$ would hold,
so that
$|\mathrm{lk}(T_{-3, 2}, c)| = |1 - \varepsilon|, |2 - \varepsilon|$, or $|3 - \varepsilon|$. 
On the other hand, 
we have $\mathrm{lk}(T_{-3, 2}, c) = 0$. 
Hence the trivial knot $c$ would be obtained
from $c_{\mu}$ by a single $(-1)$--move. 
Then, by Remark~\ref{trivializing band} the band $b$
is unique up to isotopy in $S^3$, and thus
$c = c_{\mu}\,\natural_b\,\alpha_{-1}$ is the same seiferter as 
$c^{-1}_1$ in Figure~\ref{tc1} 
($c^{-1}$ in Figure~\ref{c-1twist}). 
However, $(-1)$--twist on $T_{-3, 2}$ along $c^{-1}_1$
gives a trivial knot, 
and not $P(3, -3, -3)$.
This contradicts assertion~(2).
\QED{Lemma~\ref{P(3,-3,-3)}}

Theorem~\ref{generalizedP(3,-3,-3)} below shows that
twists along the seiferter $c$
in Figure~\ref{nonpmseiferter} extend Song's example.

\begin{theorem}
\label{generalizedP(3,-3,-3)}
Let $c$ be the seiferter for $(T_{-3, 2}, -1)$ given in Figure~\ref{nonpmseiferter}.
Then, 
all Seifert surgeries obtained from $(T_{-3, 2}, -1)$
by nontrivial twists along $c$ are 
small Seifert fibered surgeries on strongly invertible hyperbolic knots.  
However, none of them arises from
the primitive/Seifert-fibered construction.
In particular, no knot in this family has tunnel number one. 
\end{theorem}

\noindent
\textit{Proof of Theorem~\ref{generalizedP(3,-3,-3)}.}
Let $K_p$ be the knot obtained from $T_{-3, 2}$
by twisting $p$ times along $c$;
then $K_0 =T_{-3, 2}$ and $K_{-1} = P(3, -3, -3)$. 
Note that $K_p$ is the image of $T_{-3, 2}$ 
after performing $(-\frac{1}{p})$--surgery on the trivial knot $c$. 
As shown in Figure~\ref{nonpmseiferter}, 
there is a $\pi$--rotation $f$ of $S^3$ which restricts to
inversions of $T_{-3, 2}$ and $c$.
Take an $f$-invariant tubular neighborhood $N(c)$,
and extend the involution $f | S^3 - \mathrm{int}N(c)$ to 
$S^3 = S^3 - \mathrm{int}N(c) \cup_{-\frac{1}{p}} (S^1 \times D^2)$.
We obtain a strong inversion of $K_p$. 

Since $\mathrm{lk}(T_{-3, 2}, c) =0$, 
the surgery slope does not change under the twistings, 
and thus $p$--twist along $c$ converts $(T_{-3, 2}, -1)$ to
$(K_p, -1)$ for any integer $p$. 
Let $\mathcal{F}$ be a Seifert fibration of $K_0(-1)$ obtained by
extending a Seifert fibration of $S^3 -\mathrm{int}N(K_0)$
in which $s_{-3}$ is an exceptional fiber of index $3$;
$\mathcal{F}$ is a Seifert fibration over $S^2(2, 3, 5)$.
Let $\mu, \lambda$ be a preferred meridian--longitude pair of
$N(s_{-3}) \subset S^3$.
Then a regular fiber $t (\subset \partial N(s_{-3}))$ of $\mathcal{F}$
is expressed as $[t] = -3[\lambda] +2[\mu]
\in H_1( \partial N(s_{-3}) )$.

\begin{lemma}
\label{Kp-1Seifert}
$K_p(-1)$ is a small Seifert fiber space over $S^2(2, |10p +3|, 5)$.
\end{lemma}

\noindent
\textit{Proof of Lemma~\ref{Kp-1Seifert}.}
The seiferter $c_2^{-1}$ for $(K_0, -1)$ is a band sum of $s_{-3}$
and a simple closed curve $\alpha_{-1}$ in $\partial N(K_0)$
with slope $-1$ up to isotopy in $S^3 - \mathrm{int}N(K_0)$. 
Since $\mathrm{lk}(T_{-3, 2}, s_{-3}) = 2$,
by Proposition~\ref{move}(3) $0$--framing of $N(s_{-3})$ becomes
$3$--framing of $N(c_2^{-1})$ after isotoping $s_{-3}$ to $c_2^{-1}$
in $K_0(-1)$.
Since the seiferter $c$ for $(K_0, -1)$ 
is isotopic in $S^3 - \mathrm{int}N(K_0)$ to 
a band sum of $c_2^{-1}$ and
$\alpha_{-1}$, and $\mathrm{lk}(T_{-3, 2}, c_2^{-1}) =1$,
again by Proposition~\ref{move}(3)
the $3$--framing of $N(c_2^{-1})$ becomes
$4$--framing of $N(c)$ after isotoping $c_2^{-1}$
to $c$ in $K_0(-1)$. 
Hence, a preferred meridian--longitude pair of $N(c)$
 is sent to $\mu, \lambda -4\mu$ curves on $\partial N(s_{-3})$
after isotopy in $K_0(-1)$,
where $\mu, \lambda$ is a preferred meridian--longitude pair of
$N(s_{-3})$.
The $(-\frac{1}{p})$--surgery slope on $\partial N(c)$ is then
sent to $\mu -p(\lambda -4\mu) = (4p +1)\mu -p\lambda$ curve
on $\partial N(s_{-3})$.
Hence, the index of $c$ in the Seifert fibration of $K_p(-1)$
induced from $\mathcal{F}$ is 
$|( (4p +1)[\mu] -p[\lambda])\cdot (2[\mu] -3[\lambda])|
= |10 p +3|$. 
It follows that $K_p(-1)$ is a small Seifert fiber space
over $S^2(2, |10p+3|, 5)$ as desired.
\QED{Lemma~\ref{Kp-1Seifert}}

\begin{lemma}
\label{Kphyperbolic}
$K_p$ is a hyperbolic knot of genus one for $p \ne 0$.
\end{lemma}

\noindent
\textit{Proof of Lemma~\ref{Kphyperbolic}.}
The pretzel knot $P(3, -3, -3) =K_{-1}$ in Figure~\ref{tptwist}
bounds an obvious Seifert surface of genus one
disjoint from $c$.
Hence, after $(p +1)$--twist along $c$
the resulting knot $K_p$ bounds
a Seifert surface of genus one for any integer $p$.
Since Lemma~\ref{Kp-1Seifert} above implies
that $K_p$ is a nontrivial knot, $K_p$ is a knot of genus one.

If $K_p$ is not a hyperbolic knot, 
then $K_p$ is either
a satellite knot or a torus knot.
The fact that the genus of $K_p$ is one implies that
$K_p$ is a trefoil knot or a satellite knot such that
$K_p$ is null-homologous in its companion solid torus $V$.
We see from Lemma~\ref{Kp-1Seifert} that
$K_p$ is not a trefoil knot for $p \ne 0$, 
so that $K_p$ is a satellite knot.
Since $K_p(-1)$ is a small Seifert fiber space
with the trivial first homology group,
it does not contain an essential torus
\cite[Example~VI.13]{J}.
Hence, the proof of \cite[Theorem~1.4]{MM1} implies that
the manifold obtained from $V$ by $(-1)$--surgery along $K_p$
 is a solid torus.
However, this is impossible
because the winding number of $K_p$ in $V$ is zero \cite{Ga1}.
It follows that $K_p$ is a hyperbolic knot for $p \ne 0$.
\QED{Lemma~\ref{Kphyperbolic}}

\begin{figure}[htbp]
\begin{center}
\includegraphics[width=0.8\linewidth]{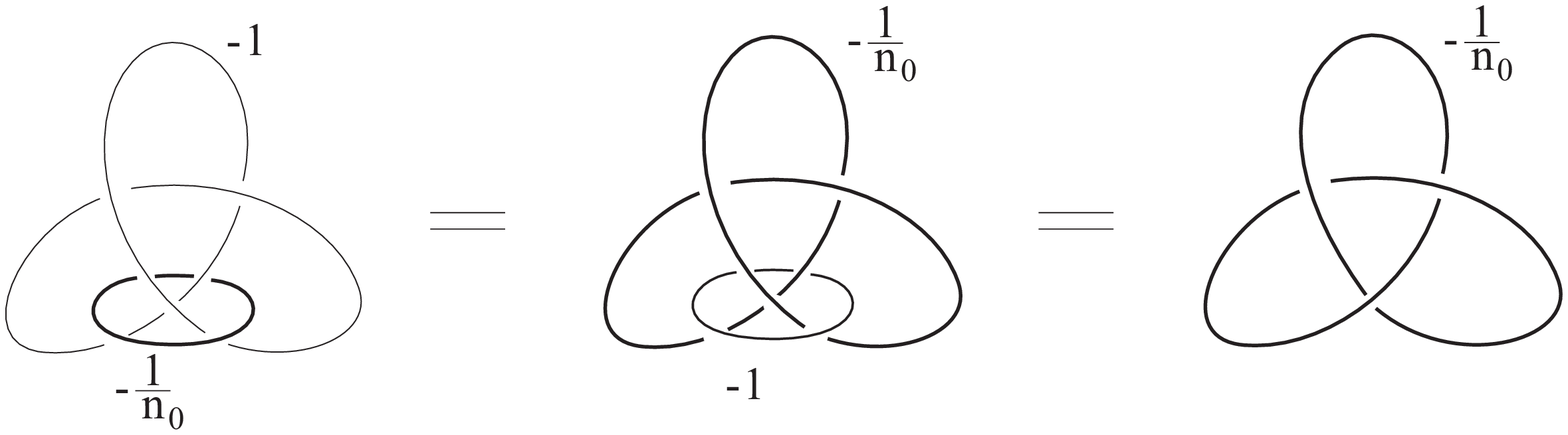}
\caption{$Tw(n_0)(-1)$ has a surgery description $T_{-3, 2}(\frac{-1}{n_0})$.}
\label{KRcalculus}
\end{center}
\end{figure}

Now suppose that
$(K_{p_0}, -1)$ arises from a primitive/Seifert-fibered
construction for some $p_0$.
Then $K_{p_0}$ has tunnel number one.
Scharlemann \cite{Sch3} has proved that knots
with both tunnel number and genus one
are 2--bridge knots or satellite knots,
as Goda and Teragaito \cite{GT} conjectured.
Since 2--bridge knots with Seifert surgeries are twist knots
\cite{BW}, $K_{p_0}$ is a twist knot or its mirror image.
In fact, $K_{p_0}$ is a twist knot $Tw(n_0)$ for some $n_0$
because $(-1)$--surgery on $K_{p_0}$ is a Seifert surgery.
Note that $(-1)$--surgery on $Tw(n_0)$ yields a Seifert fiber space
over $S^2(2, 3, |6n_0 -1|)$; 
Figure~\ref{KRcalculus} gives a pictorial proof of this fact.
On the other hand, $K_{p_0}(-1)$ is a Seifert fiber space
over $S^2(2, |10p_0 +3|, 5)$ by Lemma~\ref{Kp-1Seifert}.
It follows $|10 p_0 +3| =3$; 
then $p_0 =0$ as desired.
\QED{Theorem~\ref{generalizedP(3,-3,-3)}} 

In \cite{EJMM},
using the Montesinos trick,
we find an infinite family of small Seifert fibered surgeries on
strongly invertible hyperbolic knots which do not 
arise from the primitive/Seifert-fibered construction.
The same family is also obtained from $(T_{-3, 2}, m)$
by twisting along the annular pair $\{s_{-3}, c_1^m \}$
in Figure~\ref{T-32annularpair2}.

\bigskip

\noindent
\textbf{Acknowledgments.}
The authors would like to thank Toshiyuki Oikawa and Masakazu Teragaito for useful discussions. 
The third author would like to thank Hyun Jong Song for useful private communications, 
in particular, letting him know a remarkable example of Seifert surgery 
$(P(-3, 3, 3), 1)$. 
Finally, the authors would like to thank the referees for their careful reading and useful comments and,
in particular, for pointing out an error in the proof of Proposition~\ref{annular pairs T-32}, 
and would like to thank Mario Eudave-Mu\~noz for useful comments. 

The first and third authors were partially supported by Grants-in-Aid for JSPS Fellows (16-04787). 
The first author was 
also partially supported by JSPS Postdoctoral Fellowships for 
Foreign Researchers at Nihon University (P 04787),  
COE program ``A Base for New Developments of Mathematics into Science and Technology" at University of Tokyo, 
and INOUE Foundation for Science at Tokyo Institute of Technology.
He would like to thank 
Toshitake Kohno for his hospitality in University of Tokyo and 
Hitoshi Murakami for his hospitality in Tokyo Institute of Technology. 
The third author has been partially supported by JSPS Grants-in-Aid for Scientific 
Research (C) (No.21540098), The Ministry of Education, Culture, Sports, Science and Technology, Japan 
and Joint Research Grant of Institute of Natural Sciences at Nihon University for 2013.

\end{document}